\documentclass[11pt]{amsart}

\usepackage{amssymb} 
\usepackage[latin1]{inputenc}
\usepackage[matrix,arrow,curve]{xy}
\usepackage[mathscr]{eucal}
\usepackage{epic,eepic}
\usepackage{bbm}

\newcommand{\Q}{\mathbb{Q}}
\newcommand{\Z}{\mathbb{Z}}

\newcommand{\OO}{\mathcal{O}}

\newcommand{\G}{\mathbb G}

\newcommand{\Gl}{\mathrm{GL}}
\newcommand{\Hom}{\mathrm{Hom}}

\newcommand{\F}{\mathbb{F}}
\newcommand{\N}{\mathbb{N}}

\newcommand{\rig}{\mathrm{rig}}

\newcommand{\Ext}{\mathrm{Ext}}

\newcommand{\coker}{\text{coker}}

\newcommand{\tr}{\mathrm{tr}}

\newcommand{\rk}{\mathrm{rk}}

\newtheorem{theorem}{Theorem}
\newtheorem{prop}{Proposition}
\newtheorem{cor}{Corollary}
\newtheorem{lemma}{Lemma}
\newtheorem{question}{Question}
\newtheorem{conjecture}{Conjecture}

\theoremstyle{definition}
\newtheorem{remark}{Remark}
\newtheorem{definition}{Definition}

\def\]{\textup{\mbox{]\hspace{-.15em}]}}}
\def\[{\textup{\mbox{[\hspace{-.15em}[}}}

\newenvironment{pf}
{\medskip\noindent {\it Proof --- \ }}
{\hfill\nobreak $\Box$ \par\bigbreak}

\newcommand{\RR}{{\mathcal R}}
\newcommand{\EE}{{\mathcal E}}

\newcommand{\Ac}{{\mathcal A}}

\newcommand{\hotimes}{\hat \otimes}
\newcommand{\Fil}{\text{Fil}}
\newcommand{\gr}{\text{gr}}
\renewcommand{\sp}{\text{Sp\,}}
\newcommand{\anneau}{\mathcal O}
\newcommand{\crys}{\text{crys}}
\newcommand{\dr}{\text{dR}}
\newcommand{\ncr}{\text{ncr}}

\newcommand{\Frac}{\text{Frac}}
\newcommand{\pot}{\text{pot}}
\newcommand{\spec}{\text{Spec}\,}
\newcommand{\del}{\partial}
\newcommand{\fg}{(\phi,\Gamma)}
\newcommand{\D}{\mathcal D}
\newcommand{\Gr}{\text{gr}}
\newcommand{\Sel}{\text{Sel}}
\renewcommand{\aa}{{\underline{a}}}
\newcommand{\p}{\mathfrak{p}}

\begin{document}

\baselineskip 16pt

\title[Rank of Selmer groups in an analytic family]{Ranks of Selmer groups in an analytic family}
\author{Jo\"el Bella\"iche}
\email{jbellaic@brandeis.edu}
\address{Math Department, MS 050\\Brandeis University\\415 South Street\\ 
Waltham, MA 02453}

\begin{abstract}
We study the variation of the dimension of the Bloch-Kato Selmer group of a $p$-adic Galois representation of a number field
 that varies in a refined family. We show that, if one restricts ourselves 
to representations that are, at every place dividing $p$,
 crystalline, non-critically refined, and with a fixed number of non-negative 
Hodge-Tate weights, then the dimension of the Selmer group varies essentially 
lower-semi-continuously. This allows to prove lower bounds for Selmer groups 
``by continuity'', and in particular to deduce from a result of \cite{BC}
 that the $p$-adic Selmer group of a modular eigenform of weight $2$ of sign $-1$
has rank at least $1$.
\end{abstract}
  
\maketitle

\tableofcontents

\section{Introduction}

\par \bigskip

We first introduce some notations and conventions that shall be used throughout this paper.

Let $K$ be a number field, $p$ an odd prime number, and $\Sigma$ a finite set of finite places of $K$ containing the
set $\Sigma_p$ of places of $K$ above $p$. We denote by $G_{K,\Sigma}$ the Galois group of the maximal algebraic 
extension of $K$ unramified at all finite places that are not in $\Sigma$. For $v$ a place of $K$, we denote by $K_v$ the completion of $K$
at $v$, and by $G_{K_v}$, or $G_v$ when there is no ambiguity, 
the absolute Galois group of $K_v$. If $v \in \Sigma$, there is a natural morphism $i_v: G_v \rightarrow G_{K,\Sigma}$ 
well-defined up to conjugacy. If $G_{K,\Sigma} \rightarrow \Gl(V)$ is a representation, we will call by abuse its pre-composition with $i_v$ the {\it 
restriction} $V_{|G_{v}}$ of that representation $V$ to $G_{v}$. It is well-defined up to isomorphism.

All rings are commutative with unity. If $S$ is a ring, $\spec S$ 
is its spectrum in the sense of Grothendieck, and if $S$ is an affinoid algebra
over $\Q_p$, $\sp S$ is its attached rigid analytic space in the sense of Tate.
If $x \in \sp S$, we shall denote by $L(x)$ its field of definition, a finite extension of $\Q_p$.  
\par \bigskip

The aim of this article is to study the variation of the Bloch-Kato Selmer group
$H^1_f(G_{K,\Sigma},V_x)$ (and its closely related variants $H^1_g(G_{K,\Sigma},V_x)$ and $H^1_e(G_{K,\Sigma},V_x)$)
 when the $p$-adic representation $V_x$ of $G_{K,\Sigma}$ varies along a family of representations. 

To make this aim more precise, let us define a {\it family} of representations of a topological group $G$ of dimension $d$
 over a topological ring $S$ as a locally 
free module\footnote{We could define more generally a family
to be a continuous pseudocharacter $G \rightarrow S$ of dimension $d$, as is done in \cite{BC}, but the generality provided will not be useful here} 
$V$ over $S$ of finite rank $d$, together with a continuous $S$-linear action of $G$.
If $L$ is a topological field and $\chi : S \rightarrow L$ is a continuous and surjective morphism of algebras, $V_\chi := V \otimes_{S,\chi}
L$ is a continuous representation of $G_{K,\Sigma}$ over $L$ of dimension $d$. We shall say that the representation $V_\chi$ {\it belongs} to the family $V$.
If $S$ is an affinoid algebra, $x \in \sp S$ is a point of residue field $L(x)$ and $\chi_x : S \rightarrow L(x)$ 
is the corresponding character, we shall write $V_x$ instead of $V_{\chi_x}$.

The aim of this paper is to study how the dimensions of the spaces $H^1_f(G_{K,\Sigma},V_\chi)$ change when $\chi$ runs along a suitable set of characters
of $S$ (or of points in $\sp S$), for suitable families $V$ of representations of $G_{K,\Sigma}$, in particular those carried by eigenvarieties.

That this variation is somewhat chaotic is already exemplified by the
case of the simplest and longest studied family of Galois representations: the Iwasawa family.
Namely, take $K=\Q$ and $\Sigma=\Sigma_p$. Set $\Gamma = \Z_p^\ast$, and denote by $\Delta$ the subgroup of $(p-1)^{\text{th}}$-roots of unity in 
$\Gamma$, and by 
$\Gamma_1$ the kernel of the reduction modulo $p$ on $\Gamma$. The Teichmuller lift gives a canonical isomorphism $\Gamma=\Gamma_1 \times \Delta$.
Let $\eta_i$, for $i \in \Z/(p-1)\Z$. be the character $\Delta \rightarrow \Z_p^\ast$ that sends $t$ to $t^i$. 
 The Iwasawa algebra is $\Lambda = \Z_p[[\Gamma]]$. We have 
$\Lambda = \prod_{i \in \Z/(p-1)\Z} \Lambda_i$ with $\Lambda_i = \Lambda \otimes_{\Z_p[\Delta],\eta_i} \Z_p = \Z_p[[\Gamma_1]] \simeq \Z_p[[T]]$.
We shall work in characteristic $0$, so set $\Lambda_{\Q_p}=\Lambda \otimes \Q_p$ and $\Lambda_{i,\Q_p} = \Lambda_i \otimes \Q_p$.
If $L$ is any finite extension of $\Q_p$, a continuous character $\chi : \Gamma \rightarrow L^\ast$ defines a morphism of algebras $\Lambda_{\Q_p} \rightarrow L$ that factors through $\Lambda_{i,\Q_p}$ if and only if $\chi_{|\Delta}=\eta_i$. For $i \in \Z/(p-1)\Z$, the $i^{\text{th}}$
Iwasawa family is $V_i=\Lambda_{i,\Q_p}$ considered as a free 
$\Lambda_{i,\Q_p}$-module of rank one, with a (continuous) action of $G_{\Q,\Sigma_p}$ induced 
by the cyclotomic character $\omega_p : G_{\Q,\Sigma_p} \rightarrow \Z_p^\ast = \Gamma$. 
It is easily seen that the representation $\Q_p(n)$, that is $\Q_p$ with the action of $G_{\Q,\Sigma_p}$ given by 
$\omega_p^n : G_{K,\Sigma} \rightarrow \Q_p^\ast$, $n \in \Z$, belongs to the family $V_i$ if 
and only if $n \equiv i \pmod{p-1}$.
Assume that $i$ is odd. For $n \neq 1$ an odd integer
it is well know, as a consequence of a result of Soul\'e, that $H^1_f(G_{K,\Sigma},\Q_p(n))=H^1_e(G_{K,\Sigma},\Q_p(n))=H^1_g(G_{K,\Sigma},\Q_p(n))$ is of dimension $1$ if $n>1$ and  is $0$ if $n \leq -1$.
Therefore, there is no {\it generic values} for $H^1_f$ (or $H^1_g$ or $H^1_e$) on the families $V_i$, that is a
value that is the dimension of $H^1_f(G_{K,\Sigma},V_\chi)$ for all $V_\chi$ that belongs to $V_i$ and that are de Rham.
For the sake of completeness, let us mention that $H^1_f(\Q_p(1)) = H^1_e(\Q_p(1))=0$ while $\dim H^1_g(\Q_p(1))=1$, and also that
$\dim H^1_{\bullet}(V_\chi)=0$ (where $\bullet \in \{f,g,e\}$) when $V_\chi$ belongs to  $V_i$ and is not de Rham. 
So the ``most frequent'' dimension of $H^1_\bullet$ is $0$, 
but there are infinitely many exceptions (and moreover the set of exceptions is dense both in the Zariski topology, and the rigid analytic topology), 
a pattern which is not very usual for a family over $\Lambda_{i,\Q_p}$, a 
 Noetherian domain of dimension $1$. 

\par \bigskip

The discussion above suggests that it might be possible to prove the existence of a generic value for the dimension of the Selmer groups
of the Galois representations in a family, if one restricts our attention to the members of a family that are de Rham and that have a fixed number of
negative Hodge-Tate weights. We shall indeed be able to prove a result of this type, but only under the  additional assertion that the family is {\it 
trianguline}, and that the representations concerned are {\it non-critical}. Let us review those two notions  before stating our theorem:

\par \medskip

Let $L$ be a finite extension of $\Q_p$.
To a continuous representation of $G_{\Q_p}$ on a $L$-vector space $V$ of dimension $d$, 
one can attach (by results of Fontaine and Colmez-Cherbonnier) a
$\fg$-module $D^{\dag}(V)$ of rank $d$ (that is a free $\RR_L$-module $D$ of rank $d$ with a semi-linear action of $\phi$ and $\Gamma$, such that 
$\phi(D)$ spans $D$ over $\RR_L$, where $\RR_L$ is the Robba ring over $L$ --- see \ref{robba} for a definition) which moreover is {\it \'etale}.
 Actually $D^\dag$ is an equivalence of categories from
the category of continuous finite-dimensional representations of $G_{\Q_p}$ over $L$ to  the category of \'etale $\fg$-modules,
a full subcategory of the category of $\fg$-modules. 

According to Colmez (\cite{Co1}), a $\fg$-module $D$ over $\RR_L$ is called {\it trianguline} if it
is provided with a complete filtration $\Fil^\bullet D$ such that for $0 \leq i \leq d$, $\Fil^iD$ is a  $\fg$-submodule of $D$, 
direct summand as an $\RR_L$-module, and of rank $i$. 
Such a filtration is called a {\it triangulation}, and it may not be unique if it exists; we say that $D$ is {\it triangulable} 
to express the fact that it has at least one triangulation, but that none is specified.
 
A triangulation on $D$ defines graded pieces $\Gr^i D = \Fil^i D/ \Fil^{i-1}D$ for $i=1,\dots,d$ which are
$\fg$-modules of rank one, and Colmez has shown that those objects are in natural bijection with continuous characters $\Q_p^\ast \rightarrow L^\ast$.
A trianguline $\fg$-module thus defines a sequence of characters $\delta_i:\Q_p^\ast \rightarrow L^\ast$, $i=1,\dots,d$, that we call its {\it parameter}. 
The weight $s_i$ of the character $\delta_i$ is defined as the negative of its derivative at $1$ ({\it caveat lector}: this convention is  opposite
to the one chosen by Colmez \cite[\S2.2]{Co1}, but it is the one of \cite[\S2.3.3]{BC}). 

If $V$ is a representation of dimension $d$ over $L$  of $G_{\Q_p}$, we say that $V$ is {\it triangulable} if $D^\dag(V)$ is, and {\it trianguline} if we have 
chosen a triangulation of $D^\dag(V)$. In general, if X is an adjective describing a property of $\fg$-modules, we shall say that $V$ {\it is }X 
if and only if $D^\dag(V)$ is X.

It turns out (see \cite[Prop 2.3.3]{BC}), that if $V$ is trianguline,
the numbers $s_i$ attached to $D^\dag(V)$ are, with multiplicity, the Hodge-Tate-Sen weights of $V$ (with the standard convention that the cyclotomic character $\omega_p$ has Hodge-Tate weight $-1$). 
We shall say that a trianguline $\fg$-module is {\it non-critical} if all the $s_i$ are rational integers, and $s_1<s_2<\dots<s_d$. 
A non-critical trianguline representation is always de Rham (\cite[Prop 2.3.4]{BC}). 

An important class of examples of trianguline representations are  the crystalline representations $V$ such that the eigenvalues
of the crystalline Frobenius $\phi$ on $D_\crys(V)$ are in $L$ (the latter condition is of course 
always satisfied after replacing $L$ by a finite extension). More precisely, for such a representation $V$, any {\it refinement} of $V$, that is any full $\phi$-stable filtration of 
$D_\crys(V)$ defines a triangulation, and conversely (see \cite[Prop 2.4.1]{BC}). 
Assume for simplicity that $V$ has distinct Hodge-Tate weights. It is an easy consequences of the definitions 
that the refinement of $V$ is non-critical (in the sense of \cite[\S2.4.3]{BC}) if and only if
$D^\dag(V)$, with the attached triangulation, is non-critical (see \cite[\S2.4]{BC} for more details).

We shall say that a character $\delta : \Q_p^\ast \rightarrow L^\ast$ is {\it exceptional} if it is of the form $t \mapsto t^n$ or $t \mapsto t^n |t|$ for some integer $n \in \Z$. We shall say that a trianguline $\fg$-module $D$ is {\it non-exceptional} if all the 
characters $\delta_i$ of its parameter are not exceptional.

\par \medskip  
We have a similar notion of a {\it trianguline} family of representations $V$ of $G_{\Q_p}$ over an affinoid $S$, namely that the $\fg$-module over $\RR_S$ 
(see \S\ref{robba}) constructed by Berger and Colmez (see \cite{BeCo} and \S\ref{Ddag})  $D^\dag(V)$ has a triangulation whose  graded pieces are {\bf free} 
$\fg$-modules of rank one {\bf of character type}. Of course, if $V$ is such a family, and $x \in \sp S$ is a point of field of definition $L(x)$, then $V_x$ is a trianguline 
representation over $L(x)$. See \S\ref{trifam} for more details.
 
\begin{theorem}\label{main} Let $S$ be a an affinoid domain of dimension $1$ over $\Q_p$, and $V$  a locally free $S$-module of rank $d$ 
with a continuous $S$-linear action of $G_{K,\Sigma}$. Assume that the restriction of $V$ to $G_{K_v}$ is trianguline 
for every $v \in \Sigma_p$. For any sequence of integers $\aa = (a_v)_{v \in \Sigma_p}$, with $0 \leq a_v \leq d$ for every $v$, 
let us call $Z_{\aa}$ the set of $x \in \sp S$ such that the trianguline representation $(V_x)_{|G_v}$ is non-critical, non-exceptional, 
with exactly $a_v$ non-positive Hodge-Tate weights at every place $v$. Then there exists a natural integer  $n_\aa$, such that for every $x \in Z_\aa$, 
$\dim H^1_f(G_K,V_x) \geq n_\aa$, and for every $x \in Z_\aa$ but a finite number, $\dim H^1_f(G_K,V_x) = n_\aa$.
 \end{theorem}
 
\begin{remark}
\begin{itemize}
\item[(i)] For $x \in Z_\aa$, we have $\dim H^1_g(G_{K,\Sigma},V_x)=\dim H^1_f(G_{K,\Sigma},V_x) = \dim H^1_e(G_{K,\Sigma},V_x)$ (see Prop~\ref{pot}). So the theorem is a theorem on the $H^1_g$ and $H^1_e$ as well.

\item[(ii)]
 The hypothesis that $p$ is totally split is technical and simplify the statement. It is certainly easy to remove it, but I did not do the job. 
 Also, the hypothesis of non-exceptionality  could be weakened, and it would be possible to allow the sequence $(s_i)$, in the definition of non-critical,
to be increasing instead of strictly increasing, but this would come at the cost of a more complicated statement (separating in particular the cases of $H^1_f$, $H^1_e$ and $H^1_g$). We have chosen to present the result in the simplest formulation that seems sufficient for most applications.

\item[(iii)]

The hypothesis that $S$ has dimension $1$ could be removed 
if we knew more on the Berger-Colmez functor $D^\dag$ from family of representations over $S$ to $\fg$-modules
over $S$ (especially its commutation to arbitrary base change 
$S' \rightarrow S$, cf. \S\ref{Ddag}). See Prop~\ref{form5} to see what the 
statement of the theorem would look like for a higher dimensional $S$.

\item[(iv)]

  It may happen that for a given $\aa=(a_v)$, the set $Z_\aa$ is finite, in which case the theorem is empty.
 But when $Z_\aa$ is infinite, the theorem states that there exists a "generic value" for $\dim H^1_\bullet(G_{K,\Sigma},V_x)$
 for $x \in Z_\aa$ (and $\bullet$ being indifferently $e$, $f$ or $g$), which is obviously well-determined, and that this value is also a lower bound.

\item[(v)] It would be interesting to have some control on the set of $x \in Z_\aa$ such that $\dim H^1_\bullet (G_{K,\Sigma},V_x) > n_\aa$.
The information given by the proof does not seem useable.
 
\end{itemize}
\end{remark}

Let us just say a few words on the method of the proof. The Selmer groups are defined as the kernel of a natural map from a global
Galois cohomology group to a product of local Galois cohomology group. 
For example $H^1_f(G_{K,\Sigma},V) = \ker (H^1(G_{K,\Sigma},V) \rightarrow \prod_{v \in \Sigma_p} H^1(G_v,V \otimes B_\crys))$.
The same definition makes sense when $V$ is a family of representations over an algebra $S$ instead of an individual one. 
A formal argument, developed in \S\ref{formalism}, show that when the involved cohomology modules, and also the $H^2(G_{K,\Sigma},V)$
are finite  over $S$, the resulting Selmer groups behave well in families. Unfortunately, the modules 
$H^1(G_v,V \otimes B_\crys)$ are not finite  in general, which explains the chaotic behavior of the $H^1_f$ described at the beginning of this article.
To circumvent this problem, we use an idea due to Jay Pottharst, himself generalizing an idea of Greenberg, to replace
the $H^1_f$ by an other type of Selmer group, defined using $\fg$-modules and their cohomology, 
that we call the Pottharst Selmer group $H^1_\pot$. We prove that the cohomology groups
involved in its definition are of finite type in family (cf.~\ref{fin1} and~\ref{fin2}), and that for non-critical representations, the $H^1_\pot$ 
and $H^1_f$ are essentially the same (an easy consequence of the results of \cite[chapter 2]{BC}, and also a variant of a theorem of \cite{pot}): 
see Prop.~\ref{pot}. This, together with our formal work of \S\ref{formalism}, gives the theorem.
  
\par \bigskip

Of course, our theorem is interesting only to the extent that we are able to check that natural families are trianguline.
By natural families I mean essentially the families carried by the eigenvarieties, and more generally 
the {\it refined families} defined in \cite[Chapter IV]{BC}: they are families of representations over an affinoid algebra $S$ 
with a dense set of {\it classical points} $Z \subset  \sp S$ where the 
representations $\rho_z$ are crystalline and provided with a refinement, such that that for every $i=1,\dots,d$, the  $i^{\text{th}}$ eigenvalue 
of the crystalline Frobenius (for the order given by the refinement) 
renormalized by the $i$-th Hodge-Tate weight (in increasing order) $\kappa_i(z)$, extends to an analytic function defined on all $X$
(namely $p^{-\kappa_i(z)} F_i(z)$ extends to an analytic function over $X$). There are several other minor requirements in this definition, which is
 recalled below in \S\ref{refined}.
 
Our second main theorem (Theorem~\ref{fam2}) 
states that, if $V$ is such a family $S$ which has rank $2$, and  $x \in \sp S$ is a classical point such that $V_x$ is non-critical,
there exists an affinoid neighborhood of $x$ on which $V$ becomes trianguline. The hypothesis of non-criticality is
needed, as shown by \cite[Example 2.5.8]{BC}. 
I conjecture (cf. Conjecture \ref{conj}) that the same result holds in any rank. Note that
the formal counterpart of this statement, that is the same statement where ``an affinoid neighborhood of $x$'' is replaced by the ``completion of the
local ring of $S$ at $x$'' is known (\cite[Theorem 2.5.6]{BC}). It seems that Ruochuan Liu has recently made some important 
progresses toward this conjecture. See the discussion after it for more details.

\par \bigskip 

As an application, we can extend a theorem of \cite{BC} which states that for a modular form eigenform $k$ of weight $k=2k' \geq 4$ and of sign $-1$, 
the Selmer group $H^1_f(\Q,V_f(k'))$ has dimension at least $1$ (where $V_f$ is the $p$-adic Galois representations attached to $f$) 
to the case of forms of weight $2$, and in particular to elliptic curves\footnote{This case of elliptic curves is also a consequence of a recent result of Kim (\cite{kim}) that uses a completely different method -- I thank Emerton who pointe dout to me this reference.}.

This result, a consequence of the Bloch-Kato conjecture and of the Birch and Swinnerton-Dyer conjecture, 
was previously known only for  primes $p$ that are ordinary for $f$, with two completely different proofs: one of 
Nekovar (\cite{Nekovar}, \cite{Nekovar2}), and  one of Skinner and Urban (\cite{SU}).

\par \bigskip

 {\bf Acknowledgment.} It is a pleasure to acknowledge the influences of the work of many people on this article. Greenberg's work on 
Selmer groups in Iwasawa type families, in particular his proof of the type finiteness of Galois cohomology group in those families, has been 
influential. As it would be clear to any reader from the number of references to his work, Colmez' emphasis on $(\phi,\Gamma)$-modules
over the Robba ring,  his wonderful idea of trianguline representations, and its results about them (and their generalizations, as in the chapter 2
of my book with Chenevier, the recent work on Berger-Colmez, and various works of Ruochuan Liu) were also very important. 
But if I had to single out the most important influence on this work, it would be
Pottharst's recent work on ``triangulordinary'' representations (that is, essentially, non-critical representations in the terminology of \cite{BC}),
and especially its definition of the ``triangulordinary Selmer group''. 

I would like to thank Laurent Berger, Ga\"etan Chenevier, Pierre Colmez, Brian Conrad, Ralph Greenberg, Kiran Kedlaya and Ruochuan Liu for many useful 
conversations and emails. Also, this work was begun in May 2007 while I was invited at the University of Padua, and I would like to take this opportunity 
to thank Francesco Baldassarri and Bruno Chiarellotto for their kindness and their hospitality.

\section{Formalism of the variation of Selmer groups}

\subsection{A general formalism}

\label{formalism}
Let $R$ be a commutative domain and 
let $\Ac$ and $\Ac'$ be two abelian $R$-linear categories.
We let $(H^i_\Ac,\delta)$ be a cohomological functor from $\Ac$ to the category of $R$-modules,
and we let $(H^i_{\Ac'},\delta)$ be another cohomological functor from $\Ac'$ to the category of $R$-modules.
When no confusion can result, we shall drop the $\Ac$ or $\Ac'$ from the notations.
We say than an object $A$ of $\Ac$ is {\it torsion-free} if for every $0 \neq f \in R$,
the endomorphism $f$ of $A$ is a monomorphism. We write $A/fA$ its cokernel. Note that $f$ acts as $0$ on $A/fA$ so by functoriality, so does it on
$H^1(A/fA)$, which is therefore naturally an $R/f$-module. 
Similar definitions and considerations apply for $\Ac'$. 

Finally, we suppose given two torsion-free objects $V \in \Ac$, $V' \in \Ac'$, a morphism of $R$-modules
$u : H^1(V) \rightarrow H^1(V')$, and  for every $0 \neq f \in V$, a morphism of $R/f$-modules 
$u_f: H^1(V/fV) \rightarrow H^1(V'/fV')$  such that the 
following diagram commute 
\begin{eqnarray} \label{commute} 
\xymatrix{ H^1(V) \ar[r] \ar[d]^u &  H^1(V/fV) \ar[d]^{u_f} \\ H^1(V') 
\ar[r] & H^1(V'/fV')}
\end{eqnarray}

In this context, we define the {\it Selmer modules} $S(V)$ as $\ker u$, and $S(V/fV)$ as $\ker u_f$.
 
 From now on, we assume that $R$ is Noetherian, and call $\Pr(R)$ the set of principal
prime ideals in $R$. If $M$ is an $R$-module, we call $\rk_R M$ its {\it generic rank},
that is the dimension over $\Frac(R)$ of $M \otimes \Frac(R)$, which is a natural integer or $+ \infty$.
 
\begin{prop}\label{form1}
Let $(f) \in Pr(R)$. Assume that $H^0(V'/fV')=0$. Then the natural 
$R/f$-linear map  $S(V) \otimes_R R/f \rightarrow S(V/fV)$ is injective, and its cokernel has generic rank no more that 
$\rk_{R/f} H^2(V)[f] + \rk_{R/f} (H^1(V')/u(H^1(V)))[f]$. Moreover, if $S(V) \otimes_R R/f$, $H^2(V)[f]$ and 
 $(H^1(V')/u(H^1(V)))[f]$ are all finite  over $R/f$, then so is $S(V/fV)$.
\end{prop}
\begin{pf}
 We have the following commutative diagram where rows and columns are exact sequences.
 
 $$ \xymatrix{   & 0 \ar[d] & 0 \ar[d] & 0 \ar[d] & \\ & S(V) \ar[r]^{\times f}\ar[d] & S(V) \ar[r]^a \ar[d]& S(V/fV) \ar[d]^b & \\
  & H^1(V) \ar[r]^{\times f}\ar[d]^u & H^1(V) \ar[r] \ar[d]^u& 
 H^1(V/fV) \ar[r]^c \ar[d]^{u_f} & H^2(V)[f] \\
0\ar[r] &  H^1(V') \ar[r]^{\times f} &  
  H^1(V')\ar[r] &   H^1(V'/fV') &  }$$
 
Indeed, the commutativity of the diagram follows immediately from the commutativity of (\ref{commute}), the $R$-linearity of $u$, and the definitions of $S(V)$ and $S(V/fV)$.  The exactness of the second and third lines come from the long exact sequences of cohomology
 attached to the short exact sequence
$0 \rightarrow V \rightarrow V \rightarrow V/fV \rightarrow 0$, and similarly for $V'$ (those exact 
sequences being exact by our hypothesis that $V$ and $V'$ are torsion free), 
and the only  fact that may need a justification is the injectivity of the map labeled $\times f$
on the third line, but this comes from our hypothesis that $H^0(V'/fV')=0$.
The morphisms on the first line are defined by restriction of the morphisms of the second line, which makes the first line obviously a complex. An easy diagram-chasing argument (a part of the Snake Lemma, in fact) shows that it is exact. 

This diagram already shows that the map $S(V) \otimes_R R/f \rightarrow S(V_f)$ is injective.

From the diagram above, we deduce another one:

 $$ \xymatrix{   & 0 \ar[d] & 0 \ar[d] & 0 \ar[d] & \\ & S(V) \ar[r]^{\times f}\ar[d] & S(V) \ar[r]^a \ar[d]& \ker(c \circ b) \ar[d]^b & \\
  & H^1(V) \ar[r]^{\times f}\ar[d] & H^1(V) \ar[r] \ar[d]& 
 \ker c \ar[r]^c \ar[d] & 0 \\
0\ar[r] &  H^1(V')\ar[d] \ar[r]^{\times f} &  
  H^1(V')\ar[r]\ar[d] &   H^1(V'/fV') &  \\ & H^1(V')/u( H^1(V)) \ar[r]^{\times f} & H^1(V')/u( H^1(V)) \\ }$$
From the Snake Lemma, we deduce that the cokernel of the map $a$ (that is $\ker (c \circ b)/a(S(V))$) is isomorphic to $(H^1(V')/u(H^1(V)))[f]$. But $S(V/fV)/ \ker(c \circ b)$ is isomorphic to a sub-module of $H^2(V)[f]$.
 The proposition follows. 
 \end{pf}

\begin{remark} \label{form2} Actually the hypotheses that $R$ is a domain and that $(f)$ is prime are only used for the convenience of stating the result in term of generic ranks. Without these hypotheses, we still get that if $H^0(V'/fV')=0$ and $f$ is a monomorphism of $V$ and $V'$,
the map $S(V)/fS(V) \rightarrow S(V/fV)$ is injective, and is surjective whenever $H^2(V)[f]$ and 
$(H^1(V')/u(H^1(V))[f]$ are $0$.
\end{remark}

\begin{cor} \label{form3}
Assume that $R$ is a domain, and that $S(V)$, $H^2(V)$ and $H^1(V')$ are finite as $R$-modules.
Let $r$ be the  generic rank of $S(V)$. Then for all $(f) \in \Pr(R)$ such that $H^0(V'/fV')=0$, $S(V/fV)$ is finite over $R/f$, we have 
$\rk_{R/f} S(V/fV) \geq r$, and for almost all such $(f)$, we have $\rk_{R/f} S(V/fV) = r$.
\end{cor}
\begin{pf} By the proposition, the map $S(V) \otimes_R R/f \rightarrow S(V/fV)$ is injective for all $f$
such that $H^0(V'/fV')=0$, and by Nakayama's Lemma, $\rk_{R/f} (S(V) \otimes_R R/f) \geq \rk_R  S(V) = r$ for all $(f)$. Hence the first inequality. Also the fact that $S(V/fV)$ is of finite type follows from the proposition, since
$S(V) \otimes R/f$ is  finite over $R/f$ because $S(V)$ is finite  over $R$, and
$H^2(V)[f]$ and $H^1(V')/u(H^1(V))[f]$ are finite over $R/f$ since they are finite  over $R$ as sub-modules
of the finite  modules $H^2(V)$ and $H^1(V')/u(H^1(V))$ over the Noetherian ring $R$.
 
Since $S(V)$ is a finite $S$-module, we actually have  $\rk_{R/f} S(V) \otimes_R R/f = r$ for almost all $(f)$, and 
the finite $S$-modules $H^2(V)$ and $H^1(V')/ u( H^1(V))$ have no $f$-torsion for almost all $(f)$. The corollary follows.
\end{pf}
 
 \subsection{An important special case}
 
\label{formalism2}

In the above formalism, the morphisms $u$ and $u_f$ are {\it ad hoc}. Often, however, they appear naturally as follows.

Again $R$ is a commutative domain, and we now suppose it Noetherian, and $\Ac$ and $\Ac'$ are two abelian $R$-linear categories, 
with two cohomological functors toward the category of $R$-modules
 both denoted $(H^i)$ (or $(H^i_\Ac)$ and $(H^i_{\Ac'})$ when we really need to distinguish them).
We assume given in addition an exact $R$-linear functor $F$ from $\Ac$ to $\Ac'$, and a morphism of cohomological functor
on $\Ac$, $a: H^i_{\Ac} \rightarrow H^i_{\Ac'} \circ F$. 

Then, given an object $V$ of $\Ac$, an object $V'$ of $\Ac'$ and a morphism 
$b : F(V) \rightarrow V'$ in $\Ac'$, we can retrieve the situation of the preceding subsection, as follows. We define 
$u : H^1_\Ac(V) \rightarrow H^1_{\Ac'}(V')$ as the composition $H^1_\Ac (V) \stackrel{a}{\rightarrow} H^1_{\Ac'}(F(V)) \stackrel{H^1(b)}{\rightarrow} H^1_{\Ac'}(V') $, and $S(V)=S(V,V',b)$ as $\ker u$.

Now if $I$ is any ideal of $R$, we can define a sub-object $IV$ of $V$ as the sum in $V$ of $fV$ for $f \in I$, and therefore also a quotient object $V/IV$. Similar definitions apply to objects in $\Ac'$,
and by exactness of the functor $F$, we have $F(V/IV)=F(V)/IF(V)$.
Therefore, the map $b:F(V) \rightarrow V'$ induces a map $b_I : F(V/IV) \rightarrow V'/IV'$. By the construction of the preceding paragraph
where $R$ is replaced by $R/I$, $V$ by $V/IV$, $V'$ by $V'/IV'$ and $b$ by $b_I$, we  get a map $u_I:H^1(V/IV) \rightarrow H^1(V'/IV')$
and a $R/I$-module $S(V/IV)=\ker(u_I)$. If $J \subset I$ are ideal of $R$, we have by construction a commutative diagram:
\begin{eqnarray} \label{commute2} 
\xymatrix{ H^1(V/JV) \ar[r] \ar[d]^{u_J} &  H^1(V/IV) \ar[d]^{u_I} \\ H^1(V'/JV') 
\ar[r] & H^1(V'/IV')}
\end{eqnarray}

We shall say that a sequence $f_1,\dots,f_n$ of elements in $R$ is {\it $V$-regular} if $f_1$ is a monomorphism of $V$, $f_2$ is a monomorphism of $V/f_1V$,
 $f_3$ is a monomorphism of $V/(f_1,f_2)V$, etc. Note that when $\Ac$ is the category of $R$-modules, this agrees with the standard definition (cf. e.g. \cite{eisenbud}) of a
 regular sequence for a module $V$. In general, we shall say that an object $V \in \Ac$
 is {\it quasi-free} if every $R$-regular sequence is a $V$-regular sequence.
 Similar definitions apply for $\Ac'$.
 
\begin{prop} \label{form5} 
Assume that $V$ and $V'$ are quasi-free.
Assume that for all  ideals $I$ of $A$,
$H^i_\Ac(V/IV)$ and $H^i_{\Ac'}(V'/IV')$ are finite over $R/I$ for all $i \geq 0$ and that they are $0$ for $i > d$ for some fixed integer $d$. 

Let $r$ be the  generic rank of $S(V)$. 
Then for all prime $R$-regular ideals $\p$ of $R$ satisfying $H^0(V'/\p V')=0$, 
we have $\rk_{R/\p } S(V/\p V) \geq r$, and there exists a fixed non-zero element $g \in R$ such that equality holds if $\p $ does not contain $g$.
\end{prop}
\begin{pf} 
Let $\p$ be a prime $R$-regular ideal, and $\p=(f_1,\dots,f_n)$ where $(f_1,\dots,f_n)$ is an
 $R$-regular sequence. Let us set $I_k=(f_1,\dots,f_k)$, so $I_n=\p$.
 We first prove that if $H^0(V'/\p V') =0$, then for $0 \leq k \leq n$, 
$H^0(V'/I_kV') \otimes_R R_\p =0$. Indeed this is true for $k=n$; assume that it is true for some $k$, we shall show that it is true for $k-1$.
 We have, since $V'$ is quasi-free, an exact sequence
 $$ 0 \rightarrow V'/I_{k-1} V' \stackrel{f_k}\rightarrow V'/I_{k-1}V'
 \rightarrow V'/I_k V' \rightarrow 0.$$
 Taking the long exact sequence of cohomology, we get in particular $$H^0(V'/I_{k-1}V')  / f_k \hookrightarrow H^0(V'/I_{k-1}V').$$ Tensorizing by $R_\p$ which is flat, and using the induction hypothesis, we get that $(H^0(V'/I_{k-1} V') / f_k) \otimes R_\p =0$. Since 
 $f_k \in \p$, we conclude by Nakayama's Lemma.  
 \par \medskip
 
We now prove by induction on $k$ that for $\p$ as above, and $0 \leq k \leq n$, we have
 $$(S(V) / I_k S(V)) \otimes R_\p \hookrightarrow S(V/I_k V) \otimes R_\p.$$ For $k=0$ there is nothing to prove.
 To simplify notations, let us denote $V/I_{k-1}V$ by $\bar V$ 
 and  $V'/I_{k-1}V'$ by $\bar V'$, and let us 
 set $f=f_k$. Since $V$ and $V'$ are quasi-free, $f$ is a monomorphism of $\bar V$ and $\bar V'$. We only have to prove that $S(\bar V/f\bar V) 
\hookrightarrow S(\bar V)/f S(\bar V)$. Similarly to the proof of Prop~\ref{form1},
we get a diagram $$ \xymatrix{   & 0 \ar[d] & 0 \ar[d] & 0 \ar[d] & \\ & S(\bar V) \otimes 
R_\p
 \ar[r]^{\times f}\ar[d] & S(\bar V) \otimes R_\p \ar[r] \ar[d]& S(\bar V/f\bar V) \otimes R_\p \ar[d] & \\
  & H^1(\bar V) \otimes R_\p \ar[r]^{\times f}\ar[d] & H^1(\bar V) \otimes R_\p \ar[r] \ar[d]& 
 H^1(\bar V/f \bar V) \otimes R_\p \ar[r] \ar[d] & H^2(\bar V)[f] \otimes R_\p \\
0\ar[r] &  H^1(\bar V') \otimes R_\p \ar[r]^{\times f}  &  
  H^1(\bar V') \otimes R_\p \ar[r] &   H^1(\bar V'/f \bar V') \otimes R_\p &  }$$
using long exact sequences of cohomology tensorized by the flat $R$-module $R_\p$.
The commutativity of the diagram is obvious excepted for the commutativity 
of the lower right square, which is~(\ref{commute2}) 
and the exactness of the second line is clear,
while the exactness of the third line uses the fact that $H^0(\bar V'/f\bar V') \otimes R_\p =H^0(V'/I_kV') \otimes R_\p = 0$
that we have just proved. The exactness of the first line follows (by a part of the Snake Lemma),
so we get $(S(\bar V)/ f S(\bar V)) \otimes R_\p \hookrightarrow S(\bar V/ f \bar V) \otimes R_\p$, which completes our induction step.

Taking $k=n$, we get $(S(V) / \p S(V)) \otimes R_\p \hookrightarrow  S(V/\p V) \otimes R_\p$, but since $S(V)/\p S(V)$ and $S(V/\p V)$ are already $R/\p$-modules, we 
simply get $S(V)/\p S(V) \hookrightarrow S(V/\p V)$, so by Nakayama's lemma,
 $\rk_{R/\p} S(V/\p V) \geq \rk_R S(V)$.

\par \medskip

We now prove that this inequality is often an equality.
By the open nature of flatness, there exists a non-empty open set of $\spec R$ over which $H^i(V)$, 
$H^i(V')$ and $H^i(V')/u(H^i(V))$ are flat for all $0 \leq i \leq d$. 
Take $g$ any nonzero element in an ideal 
of definition of the closed set that is the the complement of that open set. We thus know that
 $H^i(V) \otimes_R R[1/g]$, $H^i(V')  \otimes_R R[1/g]$ and $(H^i(V')/u(H^i(V))) \otimes_R R[1/g]$ are flat over $R[1/g]$ for all $0 \leq i \leq d$, and $0$ for $i>d$.
 
Let $\p=(f_1,\dots,f_n)$ be a prime ideal, with $f_1,\dots,f_n$ a regular sequence,
 and such that $H^0(V/\p V)=0$ as above, and also assume that $g \not \in \p$.
  We shall prove by induction on $k$ the following statement :
   {\it the natural map $(S(V)/I_k S(V))  \otimes R_\p \rightarrow 
S(V/I_kV)  \otimes R_\p $ is an isomorphism, and $H^i(V/I_k V) \otimes R_\p$, $H^i(V'/I_k V') \otimes R_\p$ and $(H^i(V'/I_kV')/u_{I_k}(H^i(V/I_kV))) \otimes R_\p$ are flat over $R_\p/I_k R_\p$ for all $0 \leq i \leq d$, and $0$ for $i>d$.} 

If $k=0$, this results from the fact that $R_\p$ is a $R[1/g]$-module
and the universal nature of flatness.
 
Assume the assertion is true for $k-1$, and set as above $\bar V=V/I_kV$
(and similarly for $V'$)
and $f=f_k$. So we know that $H^i(\bar V)$, $H^i(\bar V')$ and $H^i(\bar V')/u_{I_{k-1}}(H^i(\bar V))$ are flat over $R_\p/I_{k-1}R_\p$ for all $0 \leq i \leq d$ and $0$ for $i >d$.
We prove that the same assertions hold with $\bar V$ replaced by $\bar V/ f \bar V$, by descending induction on $i$ (the result being true 
for $i>d$ by hypothesis).
We have exact sequences $$0 \rightarrow (H^i(\bar V) / f H^i(\bar V))  \otimes R_\p \rightarrow 
(H^i(\bar V/ f \bar V )) \otimes R_\p \rightarrow (H^{i+1}(\bar V)[f]) \otimes R_\p
$$
and similarly for $V'$. From the fact that $H^{i+1}(\bar V)  \otimes R_\p$ is flat over $R_\p/I_{k-1}R_\p$, in particular has no $f$-torsion, we deduce that the second map is an isomorphism, so $H^i(\bar V / f \bar V)  \otimes R_\p$ is flat over $R_\p/I_{k}R_\p$. The same holds for $V'$. By right exactness of the tensor product,  the flatness of the
module follows $$(H^i(V'/I_kV')/u_{I_k}(H^i(V/I_kV))) \otimes R_\p$$ follows from the induction hypothesis.

Finally, we can apply the proof of Proposition~\ref{form1} (and Remark~\ref{form2}) to see that 
$S(\bar V / f \bar V) \otimes R_\p =(S(\bar V)/f S(\bar V))\otimes R_\p$. 
The induction step, and therefore the proposition follows.
\end{pf}

\begin{cor} \label{form6} Same hypothesis and notations as the proposition above.
For all regular closed points $x$ of $\spec R$ such that $H^0(V'_x)=0$,
 we have $\dim S(V_x) \geq r$ with equality outside a Zariski closed subset.
 \end{cor}
 \begin{pf} 
 A closed point $x$ corresponds to a maximal ideal $\p$, and we have $V_x = V/\p V$. If $x$ is regular, the ring $R_\p$ is a regular local ring, therefore
the ideal $\p$ is $R$-regular by \cite[Cor 10.15]{eisenbud}.  
\end{pf}

\subsection{Finiteness results}
\label{fin1}

The following proposition generalizes to any $p$-complete algebra $S$ a similar result of
Greenberg for $S$ a local ring with finite residue field complete for the topology defined by its maximal ideal
(\cite[Remark 3.5.1]{G}) and an earlier result of Tate (for $S=\Z_p$ -- see \cite[\S2]{T} and \cite[Prop B.2.7]{rubin}).
A big part of its proof has been given to me by Brian Conrad. 

\begin{prop}  \label{finiteness} 
Let $G=G_{K,\Sigma}$ or $G=G_v$ for $v$ some place of $K$, or $G=I_v$ (the inertia subgroup of 
$G_v$ for $v$ a place of $K$ not dividing $p$).
 Let $S$ be a $p$-complete commutative ring, $M$ a finite free $S$-module with a continuous (for the $p$-adic topology)
 action of $G$. For any integer $i \geq 0$, let $H^i(G,M)$ denotes the continuous cohomology group of $G$  
with coefficients $M$. Then $H^i(G,M)$ is finite over $S$, and is $0$ if $i>2$. Moreover $H^i(G,M) \otimes_{\Z_p} \Q_p = H^i(G,M \otimes_{\Z_p} \Q_p)$ \end{prop}
 
\begin{pf} We claim that {\it for every finite $S$-submodule $Y$ of $H^i(G,M)$, 
$H^i(G,M)/Y$ has no non-zero $p$-divisible element}. Let $Z$ be a sub-module of finite type of $Z^i(G,M)$ (closed continuous cochains) 
that maps surjectively onto $Y$ (take for $Z$ the submodule generated by arbitrary lifts of the elements of a finite generating family 
of $Y$). Since $Z$ is finite, it is closed $p$-adically in $Z^i(G,M)$.
Now assume by that one has a family $(x_n)_{n \in \N}$ of elements in $H^1(G,M)/Y$ with $p x_{n+1} = x_n$. Choose
representatives $f_n$ of $x_n$ in $Z^i(G,M)$. We can write $p f_{n+1} = f_n + z_n + \del h_n$ with $z_i \in Z$ and $h_n \in C^{i-1}(G,M)$. From this it
follows that $f_0 = \sum_{n \geq 0} p^n z_n + \del (\sum_{n\geq 0} p^n h_n) \in Z + B^i(G,M)$ 
which shows that $x_0=0 \in H^i(G,M)/Y$.
 
Exactly as in Tate \cite[Corollary to Prop. 2.1]{T}, we can now prove the implication {\it if $H^i(G,M)/p$ is finite over $S/p$ then $H^i(G,M)$ is finite
 over $S$}. Indeed, lifting each element of a finite generating family of
$H^i(G,M)/p$ gives a finite family in $H^1(G,M)$ that generates a $S$-submodule $Y$ such that $H^i(G,M) = p H^i(G,M) + Y$. This implies that $H^i(G,M) / Y$ is $p$-divisible, hence $0$ by the claim already proved. So $H^i(G,M)=Y$ is finitely generated.

Now we prove that $H^i(G,M)/p$ is  finite  over $S/p$. Since this module is a sub-module of $H^1(G,M/p)$
we are reduced to prove that the latter is finite over $S/p$. But $M/p$ is a finite free module over $S/p$.
Remember that $G$ acts on it continuously for its quotient topology, that is the discrete topology.
 Let $e_i$ be a basis of $M/p$, and $U_i$ be a normal open subgroup of $G$ that fixes $e_i$. Then $U=\cap U_i$
 is an open normal subgroup of $G$ that acts trivially on $M/p$. By the Hochschild-Serre spectral sequence, since
 $G/U$ is finite, we are reduced to prove that the modules $H^i(U,M/p)$ are finite over $S/p$. But since $U$ acts trivially,
 we have $H^i(U,M/p)=H^i(U,\F_p) \otimes_{\F_p} M/p$, and since $H^i(U,\F_p)=\Hom_{\text{cont}}(U,\F_p)$ is well known to be finite
 (see \cite[Prop B.2.7]{rubin}), it follows that $H^i(U,M/p)$ is finite over $S/p$.
 
 So we have proved that $H^1(G,M)$ is finite over $S$. 
The fact that $H^i(G,M)=0$ if $i>2$ results, in view of the above proof, from the well-known fact that $H^1(U,\F_p)=0$ for $i>2$, where $U$ is as above (remember our running assumption that $p >2$).

 The ``moreover'' is then proved exactly like \cite[Prop. 2.3]{T}.
 \end{pf}
 
 \begin{cor} \label{fincor} Let $S$ be a reduced affinoid algebra, $V$ a locally-free finite type module over $S$ with
 a continuous action of $G$ (same notations as in the proposition). Then $H^i(G,V)$ is a finite type $S$-module, and is $0$ if $i>2$.
 \end{cor}
 \begin{pf}
 According to Chenevier's lemma (cf. \cite[Lemma 3.18]{Ch3}), there is an admissible covering of $\sp S$ by affinoid subdomains $U_n = \sp S_n$ such that
 on the subring $S_n^0$ of elements of bounded powers of $S_n$, there exists a finite free module $M_n$ with a continuous action of $G$, such that 
$M_n \otimes_{S_n^0} S_n = V \otimes_S S_n$. By the proposition, we know that
 $H^i(G,M_n)$ is of finite type over $S_n^0$, and that $H^i(G,V \otimes S_n)=H^i(G,M_n\otimes_{S_n^0} S_n)=H^i(G, M_n \otimes_{\Z_p} \Q_p)$ 
is of finite type over $S_{n_0} \otimes_{\Z_p} \Q_p = S_n$. But since $S_n$ is $S$-flat, $H^i(G,V \otimes_S {S_n}) = H^i(G,M) \otimes_S S_n$.
The $S$-module $H^i(G,M)$ is thus finite over all open subdomains of an admissible covering. By Kiehl's Theorem (\cite[Theorem 9.4.3/3]{BGR}),
it follows that $H^i(G,M)$ is finite over $S$. 
 \end{pf}

\subsection{Examples}

\subsubsection{The category $\Ac$}

\label{catAc}

In all our examples, the ring $R$ will be denoted $S$ and be either of the form $S_0 \otimes_{\Z_p} {\Q_p}$ where $S_0$ is 
a Noetherian $p$-complete algebra, or a reduced affinoid algebra over $\Q_p$. 

The category $\Ac$ is the category of finite $S$-modules, endowed with continuous $S$-linear action of $G_{K,\Sigma}$ (this makes sense since finite-$S$ module have a canonical topology, the $p$-adic topology in the first case, and the natural topology (see e.g. \cite[Prop 3.7.3/3]{BGR}) in the second case.
Our functor $H^i_\Ac(-)$ from $\Ac$ to the category of $S$-modules is simply the continuous group cohomology $H^i(G_{K,\Sigma},-)$ defined by continuous 
cochains. Note that $(H^i)_{i \in \N}$ is a cohomological functor. Indeed, for any $0 \rightarrow A \rightarrow B \rightarrow C \rightarrow 0$, the map $B \rightarrow C$ admits a continuous set-theoretical section (even a continuous $\Q_p$-linear section. In
the case where $S$ is affinoid over $\Q_p$, then $B$ and $C$ are complete $\Q_p$-vector spaces of countable type since so is $S$, and the existence of a continuous section results from~\cite[Prop. 7.2.1/4]{BGR}. The other case is an exercise),  so by standard results there is a long exact sequence of  cohomology attached to this sort exact sequence.

In particular, a family $V$ of representations of $G_{K,\Sigma}$ over $S$ is an object of $\Ac$. We fix such a family $V$.
Obviously, since $V$ is locally free as an $S$-module
it is torsion-free in $\Ac$ in the sense of \S\ref{formalism}, and even quasi-free in $\Ac$ in the sense of \S\ref{formalism2}. 

 \subsubsection{The $H^1_T$'s}
 
\label{h1t}

Let $S$ be a reduced affinoid algebra.
 Let $T$ be a subset of $\Sigma_p$. We let $\Ac'$ be the product of the 
categories $\Ac'_v$ for $v$ in $T$, where $\Ac'_v$ is the 
category of finite $S$-modules with a continuous action of $G_v$.  

Let $H^i_{\Ac'}$ be the functor of continuous local Galois cohomology: $(A'_v)_{v \in T} \rightarrow \prod_{v \in T} H^i(G_v,A'_v)$. For the same reason as for $\Ac$, the sequence $(H^i_{\Ac'})$ is a cohomological functor. 
Let $F:\Ac \rightarrow \Ac'$ be the restriction functor, which obviously is exact.
 
We set $V'=F(V)$. It is clear that $V'$ is quasi-free in $\Ac'$ in the sense of~\ref{formalism2}. In this situation, the modules $S(V/IV)$ we get (for all ideal $I$ of $S$) 
 are best noted $H^1_T(G_{K,\Sigma},V/IV) = \ker (H^i(G_{K,\Sigma},V/IV) 
\rightarrow \prod_{v \in T} H^i(G_v,V/IV))$. Two special
cases are $T = \emptyset$, where we simply get the full Galois cohomology group $H^1(G_{K,\Sigma},V/IV)$, and $T = \Sigma_p$.

Since the $H^i(G_{K,\Sigma},V)$ and $H^i(G_v,V)$ are of finite-type (and $0$ if $i >2$) by Corollary~\ref{fincor}, 
we can apply  Cor.~\ref{form6}, and we get in particular that there is an $r \geq 0$ such that  for $x$ running 
among regular closed points of $\spec S$ such that $V_x^{G_v}=0$ for $v \in T$, 
the dimension of $H^1_T(G_{K,\Sigma},V_x)$ is at least $r$
with equality for $x$ outside a Zariski closed proper subset.

\begin{remark}
This nice behavior is consistent with Jannsen's conjecture (see \cite{jannsen})
in the cases where $S$ is an affinoid over $\Q_p$ and  $T$ is either $\emptyset$ or $\Sigma_p$.  
Indeed, that conjectures says that the dimension of $H^1_T(G_{K,\Sigma},W)$, for $W$ a geometric representation of motivic weight not equal to $-1$ 
over a finite extension of $\Q_p$, and $T$ either $\emptyset$ or $\Sigma_p$, only depends on simply computable local terms and global $H^0$, 
which are easily seen to be constant on a connected family.

In particular, looking for example at a refined family (see~\ref{refined}) of representations $V$ of $G_{K,\Sigma}$ over $S$, there should be a dense
set of points $x$ such that $V_x$ is geometric (either according to the Fontaine-Mazur conjecture, or by construction of $V$ if
$V$ is the family carried by an eigenvariety), and the condition of being of weight $-1$ defines an hypersurface of $S$, since the motivic weight can be written as a linear combination of the Hodge-Tate weights at places dividing $p$. Therefore, Jannsen's conjecture predicts that the dimension 
of $H^1_T(G_{K,\Sigma},V_x)$ for $x$ geometric outside an hypersurface is constant, a prediction which is proved by the above result, and extended to non-geometric points. 

In the intermediate cases where $T$ is not empty nor $\Sigma_p$, I have heard of no conjecture predicting the dimension of $H^1_T(G_{K,\Sigma},W)$,
but the above results suggest that such a conjecture, generalizing Jannsen's conjecture, should exist and have a simple form.

I plan to come back to these questions in a subsequent work.
 \end{remark}

\subsubsection{The $H^1_f$, $H^1_e$, $H^1_g$'s}

\label{h1f}
Our second example is in fact a set of three counter-examples, that we will not develop in full details.

Let $\Ac'$ be the product over $v$ in $\Sigma_p$ of categories $\Ac'_v$, where $\Ac'_v$ is a suitable abelian category of 
topological modules over $S$ with a continuous $G_v$-action. Let $H^i$ be the continuous cohomology functor as in the above example, 
  and let $F_f:\Ac \rightarrow \Ac'$ (resp. $F_e$, $F_g)$ be the functor 
  $V \rightarrow (V_{|G_v} \otimes B_\crys)_{v \in \Sigma_p}$ (resp., with $B_\crys$ replaced by $B_\crys^{\phi=1}$, $B_{\dr}$)
  which obviously is exact. 
  
  In this case, we write $H^1_f(G_{K,\Sigma},V/IV)$ (resp. $H^1_e$, resp. $H^1_g$) instead of
 $S(V/IV)$.
 
 In those cases, the conclusion of Corollary~\ref{form3} in general does not hold. The simplest way to see it
 is for $H^1_e$ for the family $\Lambda_{i,\Q_p}$ (with $i$ odd) described in the introduction. In this case, $K=\Q$, and $\Sigma=\Sigma_p$. As we have seen,
 for every $n \equiv i \pmod{p-1}$, the representations $\Q_p(n)$ belong to that family.
But for those representation, we have $H^1_e(G_{\Q,\Sigma_p},\Q_p(n))=0$ if $n < 0$ and $\dim H^1_e(G_{\Q,\Sigma_p},\Q_p(n))=1$ 
if $n$ is positive, contradicting the conclusion of Corollary~\ref{form3}, since obviously, we have $H^0(G_p,\Q_p(n) 
\otimes B_\crys^{\phi=1})=D_\crys(\Q_p(n))^{\phi=1}=0$ (for $n \neq 0$).
 One can infer that some hypothesis of the corollary does not hold. Since $H^i_{\Ac}(\Lambda_{i,\Q_p})=H^i(G_{\Q,\Sigma_p},\Lambda_i \otimes_{\Z_p} \Q_p)$ 
is finite for all $i$ 
by Prop~\ref{finiteness}, we easily see that the hypothesis that is not satisfied here is the one that asks
 that $H^1(G_v, V \otimes B_\crys^{\phi=1})$ be of finite type over $\Lambda_{i,\Q_p}$.
 
\subsubsection{Greenberg Selmer groups}

Let $S$ be an affinoid algebra over $\Q_p$.
In this example, we let $\Ac'$ be the product of the categories $\Ac'_v$ for $v$ in $\Sigma_p$, where $\Ac'_v$ is the 
category of finite $S$-modules with a continuous action of $G_v$, and the functor $H^i_{\Ac'}$ and $F$ are as in  \S\ref{h1t}.
But instead of taking $V'=F(V)$, we start with a $V'$ which is a quotient of $F(V)$. More precisely, we assume that $V'=(V'_v)_{v \in \Sigma_p}$ 
where each $V'_v$ is a free $S$-module with a continuous action of $G_v$ which is a quotient of $V_{|G_v}$. If we call $b$ the canonical morphism $F(V) \rightarrow V'$ in $\Ac'$, we are exactly in the situation of~\ref{formalism2}. We can thus define for every ideal $I$ of $S$, 
Selmer groups $S(V/IV)$ that we can call the {\it Greenberg Selmer Group} $H^1_{\gr}(G_{K,\Sigma}(V/IV))$ of $V/IV$ and $(V'/IV')$: we have
$S(V/IV)= H^1_{\gr,V'}(G_{K,\Sigma},V) = \ker (H^1(G_{K,\Sigma},V) \rightarrow \prod_{v \in \Sigma_p} H^1(G_v,V'_v))$.

 The finiteness hypotheses of Prop.~\ref{form5} are satisfied by Corollary~\ref{fincor}. Therefore, we see that for all points $x$ 
in $\sp S$ outside a Zariski closed proper subset, the dimension of $H^1_{\gr,V'}(G_{K,\Sigma},V_x)$ is constant. A remarkable fact, due to Flach,
is that for suitable choices of $V'$ and $x$ (namely such that the Hodge-Tate weights of $(V'_v)_x$ are exactly the Hodge-Tate weights of $(V_v)_x$) 
then we have $H^1_g(G_{K,\Sigma},V_x) = H^1_{\gr,V'}(G_{K,\Sigma},V_x)$. This allows to get some information about the behavior of the $H^1_g$ not given by the direct approach~\ref{h1f},
provided we can construct a suitable quotient $V'$ of $F(V)$. Unfortunately, the conditions under which we can construct such a $V'$, 
known as the Panchiskin condition, are too restrictive for the applications to eigenvarieties (except on their ordinary loci), 
and the aim of this paper is precisely to remove them.   

 \section{Trianguline families of $\fg$-modules}

\subsection{Families of $\fg$-modules and their triangulations}

\label{trianguline}

\subsubsection{The relative Robba ring}

\label{robba} 

Let $S$ be any ring with a non archimedean complete valuation $|\ |$. For any nonnegative real $\rho < 1$, we set 
$$\RR_S^\rho = \{ f = \sum_{i=-\infty}^{+\infty} a_i T^i \, | \, a_i \in S, 
f(T)\text{ is convergent for }\rho \leq |T| < 1\}.$$ 
The convergence condition means that $|a_i| \rho^i$ goes to $0$ when $i$ goes to $-\infty$ and that for every positive real number $r<1$, $a_i r^i$ goes to 
$0$ when $i$ goes to $+\infty$. Therefore, for any real $r$ such that $\rho < r <1$, and any $f \in \RR_S^\rho$, we can set 
$|f|^{(r)} = \sup_{i \in \Z} |a_i r^i|$.

For $f = \sum_{i \in \Z} a_i T^i$ and $g=\sum_{\i \in \Z} b_i T^i$ two elements of 
$\RR_S^{\rho}$, then for every $i$ the sum 
$c_i = \sum_{j+k = i} a_j b_k$ converges, and the series $\sum c_i T^i$ is an element of $\RR_S^\rho$ that we call $fg$.
\begin{lemma} The obvious addition and the multiplication just defined makes $\RR_S^\rho$ a commutative ring, and the $|\ |^{(r)}$ are valuation on that ring.
If $S$ is a Banach algebra over $\Q_p$, then those valuations 
make $\RR_S^{\rho}$ a Frechet space over $\Q_p$, and a topological ring for its
Frechet topology.
\end{lemma}

Note that for $1 > \rho' > \rho$, we have a continuous inclusion $\RR_S^\rho \rightarrow \RR_{S}^{\rho'}$.
We set $$\RR_S = \bigcup_{\rho <1} \RR_S^{\rho}.$$
The set $\RR_S$ is obviously a ring (and an $S$ algebra), and it is called the {\it (relative) Robba ring} over $S$. When $S$ is a Banach algebra over 
$\Q_p$, $\RR_S$ has a natural topology as an inductive limit of the Frechet space $\RR_S^\rho$.

\begin{lemma} \label{liu} 
(cf. Liu \cite[Prop 3.2]{Liu1}) If $S$ is a reduced affinoid algebra over $K$ (with its spectral norm, which is a complete valuation), 
then $\RR_S^\rho$ (resp. $\RR_S$) 
is naturally isomorphic as a topological $S$-algebra to $\RR_{\Q_p}^\rho \hotimes S$
(resp. to $\RR_{\Q_p} \hotimes S$).
\end{lemma}

We shall also deal with  the following subrings of $\RR_S$:
$\EE^\dag S:=\{f=\sum a_i T^i \in \RR_S,\ (|a_i|)\text{ is bounded}\}$, $\RR^+_S := \{s = \sum_i a_i T^i \in \RR_S,\ a_i =0 \text{  for }i<0\}$,
and $\EE^+_S := \EE^\dag_S \cap \RR_S^+$.

\subsubsection{Different notions of families of $\fg$-modules}
\label{different}

Let $S$ be an affinoid algebra over $\Q_p$.

By a {\it general $\fg$-module} over $S$ we mean a module $D$ over $\RR_S$ with a semi-linear action of $(\phi,\Gamma)$.
Note the absence of conditions of finiteness or freeness of $D$, or of a topology on $D$. The categories of {\it general $\fg$-modules} over $S$,
is an abelian $S$-linear category (but not $\RR_S$-linear). 

By a {\it family of $\fg$-modules} (of rank $d$) over $S$, we mean a general $\fg$-module $D$ over $S$, such that $D$ is  finite and projective
over $\RR_S$ (or, which amounts to the same, finite presentation and flat, or finite presentation and locally free), such that for its natural topology
(as a sub-module of some free finite module over $\RR_S$), the action of $\phi$ and $\Gamma$ are continuous, and such that $\phi(D)$ spans $D$ over $\RR_S$.
We say that a family $D$ of $\fg$-modules is {\it free} if $D$ is free over $\RR_S$.  

If $f : S \rightarrow S'$ is a morphism of reduced affinoid algebras, and $D$ is a general $\fg$-module over $S$, then $D_{S'}:= D \otimes_{\RR_S} \RR_{S'}$ has an obvious structure of general $\fg$-module over $S'$. If $D$ is a family (resp. a free family), then so is $D_{S'}$. We often adopt a geometric notation: If $X'=\sp S'$ and $X=\sp S$, we write $D_{X'}$ instead of $D_{S'}$ and we call it the pull-back of $D$ by the morphism $X' \rightarrow X$ induced by $f$. 
We note that If $D$ is a free family of $\fg$-modules over $S$, then $D$ has a natural structure of Ind-Frechet vector pace over $\Q_p$, so we can talk about 
$D \hotimes_S S'$ and we have, using Lemma~\ref{liu}, 
\begin{eqnarray} \label{DS'} D_{S'} = D \hotimes_S S'. \end{eqnarray}

The theory of families of $\fg$-modules would be much simplified if we had an answer to the following questions (especially if we had a positive answer to the first one):

\begin{question}\label{Q1} If $D$ is a family of $\fg$-modules over $S$ of rank $d$, does there necessarily exist a covering of $S$ by affinoid 
algebras $S_n$ such that the families of $\fg$-modules $D_{S_n}$ are free? If so, can we choose the covering ($S_n$) admissible?
\end{question}

\begin{remark}
We recall that if $S$ is an affinoid algebra over $\Q_p$ that is a field, the answer to this question is yes. Indeed such a field is a finite extension $L$  
of $\Q_p$, and the freeness of any finite projective module over $\RR_L$ is a consequence of its being a Bezout domain, which results from the first part of Lazard's paper
\cite{lazard} (the one concerned with a discrete valuation, and accessible by what Lazard calls "multiplicative method").
 
As Kiran Kedlaya made me notice, we can ask the same question for $S$ that are Banach field extensions of $\Q_p$ but not finite 
(so not affinoid algebra over $\Q_p$), and in its generality the answer may be false. Indeed,
the hard part of \cite{lazard} (with non-discrete valuation) shows that the answer is yes if and only if the field $S$ is spherically complete. This may be an indication that the question is difficult. 
On the other hand, the spectral valuation on a reduced  affinoid algebra is discrete, which may help adapt the proof from the ``easy'' part of \cite{lazard}.
\end{remark}

If $D$ is a family of $\fg$-modules over $S$, and $x \in \sp S$, we denote by $D_x$ the $\fg$-module $D_{L(x)}$ over $
\RR_{L(x)}$ equal to  $D \otimes_{\RR_S} \RR_{L(x)} = D \otimes_S \otimes L(x)$. It is  a projective, hence free
 $\RR_{L(x)}$-module.

\subsubsection{The functor $D^\dag$}

\label{Ddag}

Let $S$ be a reduced affinoid algebra. A construction of Berger and Colmez \cite{BeCo} attaches to any family of representations $V$ of $G_{\Q_p}$ over $S$
a family of $\fg$-module $D^\dag(V)$, at least when $V$ admits a model $V_0$ over the algebra $S_0$ of power-bounded elements of $S$, that is a finite locally free $S_0$-module with a Galois action such that $V_0 \otimes_{S_0} S = V$. Chenevier
has shown (\cite[Lemma 3.18]{Ch3}) using Raynaud's technics  that, at least locally (that is on the constituents on an admissible affinoid covering of $\sp S$) there always exists
such a model $V_0$. This allows to construct $D^\dag(V)$ in general by gluing, and
it is easy to prove that the result is independent of the choices made.
Moreover, if $f : V \rightarrow V'$ is a morphism of representations over $S$,
then locally one can chose models $V_0$ and $V'_0$ over $S_0$
and a morphism $f_0:V_0 \rightarrow V'_0$ such that $f_0 \otimes S = f$. This results from the same Raynaud's technics used by Chenevier. From that we deduce that
$D^\dag$ is a {\it functor} from the category of families of 
representations over $S$ to the  category of families of $\fg$-module over $S$
(the fact that it is a functor is stated without proof in \cite{BeCo}). 

Liu (cf. \cite{Liu1}) has shown that this functor is fully faithful and it is clear that it transforms exact sequences into exact sequences.

\begin{question} \label{QDdag} Can the functor $D^\dag$ be extended as an exact functor from the category $\Ac$ (see\S\ref{catAc}) of finite $S$-modules
with a continuous $G_{\Q_p}$-action to the category of general $\fg$-modules?
\end{question}
A  simpler question is:
\begin{question} \label{QDdag2} If $\pi : S \rightarrow S'$ is a surjective morphism of
affinoid algebras, and $V$ a family of representations over $S$, do we have a canonical
isomorphism
$(D^\dag(V))_{S'} \simeq D^\dag(V \otimes S')$?
\end{question}
A positive answer to question~\ref{QDdag} implies a positive answer to 
question~\ref{QDdag2}: let $\ker \pi=(f_1,\dots,f_n)$, and consider the $\Ac$-exact sequence
 $ V^n \rightarrow V \rightarrow V/IV \rightarrow 0$ where the first map is given by $ (v_1,\dots,v_n) \mapsto f_1v_1 + \dots + f_n v_n$. If $D^\dag$ is an exact functor from $\Ac$, we get an exact sequence of $\fg$-modules over $S$ : $D^\dag (V)^n \rightarrow
 D^\dag(V) \rightarrow D^\dag(V/IV) \rightarrow 0$ hence a canonical isomorphism 
 $D(V \otimes S') = D(V/IV) \simeq D^\dag(V)/I D^\dag(V) = D^\dag(V) \otimes_S S'
 \simeq D^\dag(V)_{S'}$, the last equality being~(\ref{DS'}).
 
\subsubsection{Free families of $\fg$-modules of rank one}

\label{free1}

Let $\delta : \Q_p^\ast \rightarrow S^\ast$ be a continuous character.
We define a free family of  $\fg$-modules $\RR_S(\delta)$ of rank one as being $\RR_S$ as a $\RR_S$-module with semi-linear action of $\phi$ and 
$\Gamma$ given by $\phi(1):= \delta(p),\  \gamma(1):=\delta(\gamma), \forall \gamma \in \Gamma$

We shall say that the family of $\fg$-module $\RR_S(\delta)$ is {\it of character type}.

\begin{lemma}~\label{rankone} Let $D$ be a free family of $\fg$-modules of rank one over $S$, and $x \in \sp S$. Then there exists an affinoid domain $\sp S' \subset \sp S$ 
containing $x$ such that the restriction $D_{S'}$ of $D$ to $S'$ is isomorphic
to $\RR_{S'}(\delta_{S'})$ for a unique continuous character $\delta_{S'}:\Q_p^\ast \rightarrow S'^\ast$ 
\end{lemma}
\begin{pf} By twisting by a constant character 
we can reduce to the case where $D_x$ is etale. Then a theorem of Liu (\cite{Liu2})
shows that there exists  an affinoid domain $\sp S' \subset \sp S$ 
containing $x$ such that the restriction $D_{S'}$ of $D$ to $S'$ is isomorphic
to $D^\dag$ of a continuous character of $G_{\Q_p}$ to $S'^\ast$. By local class field theory, this characters defines a $\delta_{S'}:\Q_p^\ast \rightarrow S'^\ast$ and 
it is clear that $D_{S'}=\RR_{S'}(\delta_{S'})$.
\end{pf}

\begin{question} \label{Q2} 
If $D$ is a  free family of $\fg$-modules of rank one over $S$,
is it globally of character type?
\end{question}
Of course, it is locally in a neighborhood of every point $x$, but there is a priori 
no reason that using those neighborhoods one could make an admissible covering, which would allow to glue the characters $\delta_{S'}$.

\subsubsection{Trianguline families}

\label{trifam} 
By a {\it  triangulation} of a family $D$ of $\fg$-modules over $S$, we mean an increasing filtration $(\Fil_i D)_{i \in \N}$ of  
subfamilies, such that $\Fil^i D$ is a family of rank $i$ for $i=0,\dots,d$ 
and such 
that for all $i=1,\dots,d$, $\gr^iD:=Fil^i D / \Fil^{i-1} D$ is of character type.

A family of  $\fg$-module provided with a triangulation (there may exist several) is called a  {\it trianguline} family.  A family of $\fg$-modules that admits a triangulations is called a {\it triangulable}. Obviously, a triangulable family is free.

\subsection{Triangulation of some refined families}

\subsubsection{Definition of a refined family}

\label{refined}
We take the following definition from~\cite[Chapter 4]{BC} except that instead of a 
pseudocharacter, we assume that our family is given by a free module with Galois
action. The greater generality given by a pseudocharacter would be of no use here.

\begin{definition}\label{defrefi} Let $X=\sp S$ be a reduced  affinoid space over $\Q_p$,
A {\it (rigid analytic) family of refined $p$-adic representations} (or shortly, a {\it refined family}) 
of dimension $d$ over $X$
is a locally free module $V$ of dimension $d$ over $S$, together with a continuous $S$-linear action of $G_{\Q_p}$ and with the following data
\begin{itemize}
\item[(a)] $d$ analytic functions $\kappa_1,\dots,\kappa_d \in \anneau(X)=S$,
\item[(b)] $d$ analytic functions $F_1,\dots,F_d \in \anneau(X)=S$,
\item[(c)] a Zariski dense subset $Z$ of $X$;
\end{itemize}
subject the following requirements.
\begin{itemize}
\item[(i)] For every $x \in X$, the Hodge-Tate-Sen weights of $V_x$ are, with multiplicity, $\kappa_1(x),\dots,\kappa_d(x)$.
\item[(ii)] If $z \in Z$, $V_z$ is crystalline (hence its weights $\kappa_1(z),\dots,\kappa_d(z)$ are integers).
\item[(iii)] If $z \in Z$, then $\kappa_1(z)<\kappa_2(z)<\dots<\kappa_d(z)$.
\item[(iv)]  The eigenvalues of the crystalline Frobenius acting on $D_\crys(V_z)$ are distinct and are
$(p^{\kappa_1}(z) F_1(z),\dots, p^{\kappa_d}(z) F_d(z))$.
\item[(v)]  For $C$ a non-negative integer, let $Z_C$ be the set of $z \in Z$ such that 
$$\kappa_{n+1}(z) - \kappa_n(z) > C \,(\kappa_{n}(z)-\kappa_{n-1}(z)) \,\, \, \text{for all } \,n\,=\,2,\dots,d-1,$$
and $\kappa_{2}(z)-\kappa_1(z) > C$. Then for all $C$, $Z_C$ accumulates at any point of $Z$. In other words, for all $z \in Z$ and $C>0$, there is a basis of 
affinoid neighborhoods $U$ of $z$ such that $U 
\cap Z$ is Zariski-dense in $U$.

\item[($\ast$)] For each $n$, there exists a continuous character $\Z_p^* 
\longrightarrow \OO(X)^*$ whose derivative at $1$ is the map $\kappa_n$ and whose evaluation at any point $z \in Z$ is the elevation to the $\kappa_n(z)$-th power.
\end{itemize}
\end{definition}

Recall that at every $z \in Z$, the representation $V_z$ is crystalline, and
has a natural refinement (that is an ordering of its crystalline Frobenius eigenvalues),
namely $(p^{\kappa_1}(z) F_1(z),\dots, p^{\kappa_d}(z) F_d(z))$. We call this refinement the {\it canonical} refinement of $V_z$.

We call $Z_\ncr$ the set of $z \in Z$ such that the canonical refinement
of $V_z$ is non-critical (see \cite[\S2.4.3]{BC}),
that is the filtration on $D_\crys(V_z)$ given by the refinement is in generic position with respect to the Hodge-Tate filtration. By condition (v) and \cite[Remark 2.4.6(ii)]{BC}, one sees easily that the set $Z_\ncr$ accumulates at any point of $Z$.

Also note that if $x \in Z$, and $U \subset X$ is any affinoid subdomain of $X$, then the restriction of $V$ to $U$ (together with the restricted data $Z \cap U$, $(F_i)_{|U}$, etc.) is again a refined family over $U$. To make this true is indeed the main reason
of the complicated form of condition (v), instead of simply asking that $Z_C$ is Zariski-dense. We shall use this remark freely in the sequel.

\subsubsection{Triangulation of  families of rank 2}

\begin{theorem}\label{fam2} Let $(X,V)$ be a refined family of representations of $G_{\Q_p}$ of dimension $2$. 
Let $x \in Z_\ncr$ such that $V_x$ is irreducible, and such that $F_1(x)/F_2(x) \not \in 
p^\Z$.
Then there is an affinoid subdomain $U$ of $X$ containing $x$ such that 
$D^\dag(V_{U})$ is trianguline, and the induced triangulation of $D^\dag(V_x)$ is the
one given by its canonical refinement.
\end{theorem}
\begin{pf}
For $i=1,2$, let $\delta_i:\Q_p^\ast \rightarrow S^\ast$ defined by $\delta_i(p)=F_i$ and
$(\delta_i)_{\Gamma}$ is the character $\Z_p^\ast \rightarrow S^\ast$ whose derivative at $1$ is $-\kappa_i$ and evaluation at any $z \in Z$ is $t \mapsto t^{-\kappa_i(z)}$;
The existence of such a character is property (*) above, and its uniqueness is obvious.

\begin{lemma} Up to replacing $X=\spec S$ 
by a smaller affinoid containing $x$, there exist
\begin{itemize}
\item[(i)] A triangular family $D$ of $(\phi,\Gamma)$-modules over $S$
such that for every $y \in \sp S$, $D_y$ is a non trivial extension
of $\RR_{L(y)}(\delta_{2,y})$ by $\RR_{L(y)}(\delta_{1,y})$.
\item[(ii)] A free module $W$ of rank $2$ over $S$ with a continuous action of $G_{\Q_p}$ such that $D^\dag(W) \simeq D$.
\end{itemize}
\end{lemma}
\begin{pf}
Note that by weak admissibility and irreducibility of $V_x$, we have
$\kappa_1(x) < v_p(p^{\kappa_i(x)} F_i(x)) < \kappa_2(x)$ for $i=1,2$. It follows that 
$v_p(F_2(x)) < 0 < v_p(F_1(x))$, so $v_p(\delta_{2,x}(p)) < v_p (\delta_{1,x}(p))$.
By continuity, after shrinking $X$ if necessary, we can assume that this property 
holds for all $y$ in $X$. We also have, for all $y \in X$, that $\kappa_1(y) + 
\kappa_2(y) = v_p ( p^{\kappa_1(y)} F_1(y) p^{\kappa_2(y)} F_2(y))$, so we have  $v_p(\delta_{2,y}(p)) + v_p (\delta_{1,y}(p)) =0$.

Colmez calls $\SS_\ast$ the rigid analytic moduli space of characters 
$\delta_1,\delta_2$ satisfying the condition $v_p(\delta_2(p)) < v_p(\delta_1(p))$ and 
$v_p(\delta_2(p)) + v_p(\delta_1(p))=0$. (See \cite[\S0.2]{Co1}, but note that he has weights conventions opposed to ours). That is to say, the characters $\delta_1,\delta_2$
define an analytic map $\sp S \rightarrow \SS_\ast$, such that $\delta_i$ is the pull-back of the universal character $\tilde \delta_{i}$ for $i=1,2$.

In \cite[\S5.1]{Co1}, Colmez constructs, locally on $\SS_\ast$ near any point
where $\tilde \delta_1 \tilde \delta_2^{-1} \not \in p^{\Z}$, a family of $\fg$-modules
of rank $2$ whose fiber at every point $y$ is a non-trivial extension of  
$\tilde \delta_{2,y}$ by $\tilde \delta_{1,y}$. The pull-back $D$ to $S$ of this 
families of $\fg$-modules obviously satisfies (i).

To prove (ii), we note that by construction $D_y$ is an \'etale $\fg$-module. Therefore,
by \cite[Main Theorem]{Liu2}, there exists a representation $W$ on $S$ (up to
shrinking it again) such that $D^{\rig}(W)=D$.
\end{pf}

To finish the proof of the theorem, it is enough to prove that $V \simeq W$ 
as $G_{\Q_p}$-representations over $X$, after shrinking $X$ if necessary. But 
\begin{lemma} There exists an affinoid neighborhood $U$ of $x$, such that for every $z \in Z_\ncr \cap U$, we have $V_z \simeq W_z$ as $G_{\Q_p}$-representations.
\end{lemma}
\begin{pf} By construction, $D^\dag(W_y)$ 
is a non-trivial extension of $\RR_{L(y)}(\delta_{2,y})$ by $\RR_{L(y)}(\delta_{1,y})$
for every $y \in X$. By \cite[Theorem 2.9]{Co1}, $\Ext^1(\RR(\delta_{2,y}),\RR(\delta_{1,y}))$ has dimension $1$ for $y$ in an affinoid neighborhood of $x$, and 
shrinking $X$ if necessary, for all $y$. Therefore, there exists up to isomorphism
only one $\fg$-module that is a non-trivial 
extension of  $\RR_{L(y)}(\delta_{2,y})$ by $\RR_{L(y)}(\delta_{1,y})$. Since $D^\dag$
is an equivalence of category (over a base field), we only have to show
that for $z \in Z_\ncr$, $D^\dag(V_z)$ is such an extension.

By Kisin's theorem (\cite[Corollary 5.16(c)]{kis}, or \cite[\S3.3]{BC} in this context), 
we have
$$D_\crys(V_z)^{\phi= p^{\kappa_1(z)} F_1(z)} \neq 0.$$ Therefore, by \cite[Prop 4.3]{Co1},
there exists an integer $k \in \Z$ such that if 
$\delta : \Q_p^\ast \rightarrow L(z)^\ast$ is the character such that $\delta$ is $t \mapsto t^k$ on $\Gamma=\Z_p^\ast$ and $\delta(p)=p^{k + \kappa_1(z)} F_1(z)$, then $D^\dag(V_z)$ has $\RR(\delta)$ as a saturated sub-$\fg$-module.
The integer $k$ is not explicitly determined in \cite{Co1}, but actually it is easy to
see that it is the opposite of the 
weight of the line $D_\crys(V_z)^{\phi= p^{\kappa_1(z)} F_1(z)}$ (see \cite[Prop 2.4.2]{BC}),
that is, since $V_z$ is non-critical, $k=-\kappa_1(z)$. Thus $\delta = \delta_{1,z}$, and
since $\RR(\delta_{1,z})$ is a saturated $\fg$-submodule of $D^\dag(V_z)$, 
the quotient is a $\fg$-module of rank one which, by looking at the determinant is seen
to be $\RR(\delta_{2,z})$. Hence $D^\dag(V_z)$ is an extension of $\RR(\delta_{2,z})$ by $\RR(\delta_{1,z})$ and this extension is non trivial, since otherwise, it would not be \'etale.
\end{pf}

Going back to the proof of the theorem, we see using the lemma that $\tr V_z = \tr W_z$
for all $z \in Z_\ncr \cap U$. Since (restricting $U$ again if necessary),
$Z_{\ncr} \cap U$ is Zariski-dense in $U$, we see that $\tr V = \tr W$ as maps $G_{\Q_p} \rightarrow S$. Since $\anneau_{X, y}$ (the rigid analytic
local field of $X$ at $y$) is Henselian, and since $V_y$ is irreducible,
a famous theorem of Serre and Carayol says that $V$ and $W$ are isomorphic over $\anneau_{X,y}$. Therefore, they are isomorphic over some affinoid neighborhood of $y$, which is what we wanted to prove.

Finally, the triangulation on $D_x$ that we have constructed is the unique one
that has $\RR(\delta_{1,x})$ as $\Fil^1 D_x$. Since the refinement of $V_x$ is 
non-critical, the attached triangulation of $D^\dag(V_x)$ is the one whose $\Fil^1$
is also $\RR(\delta_{1,x})$ (cf. the proof of the above lemma). This proves the last assertion of the theorem.

\end{pf}

\begin{remark} \label{remfam2}
If $V_x$ is reducible a variant of this theorem has been known for long. Actually, if $V_x$
is reducible, then either $v_p(F_1(x))$ or $v_p(F_2(x))=0$. If $v_p(F_1(x))=0$, we say that the point $x$, or the canonical refinement of $V_x$, is {\it ordinary}. In this case,
the property $v_p(F_1)=0$ continues to be true on some neighborhood $U$ of $x$, and on $U$, the character $\delta_1$ and $\delta_2$ are of Galois type (that is, $\delta_i(p)^n$
goes to $1$ when $n$ goes multiplicatively to $\infty$ -- after shrinking $U$ again if necessary) and the family $V$ is well-known to be an extension of the character $\delta_2$ by the character $\delta_1$. This obviously implies, and is much stronger than, that
$D^\rig(V)$ on $U$ is trianguline. 

Note also that in the case where $V_x$ is reducible and semi-simple, a refinement of $V_x$ is ordinary if and only if it is critical. If $V_x$ is not semi-simple, only one refinement is ordinary but both are critical. 
\end{remark}

\subsubsection{A conjecture about triangulation of refined families}

\begin{conjecture} \label{conj}
Let $(X,V)$ be a refined family of representations of $G_{\Q_p}$ of any 
dimension $d$. 
Let $x \in Z_\ncr$. Then there is an affinoid subdomain 
$U$ of $X$ containing $x$ such that  $D^\dag(V_{U})$ is trianguline, 
and the induced triangulation of $D^\dag(V_x)$ is the
one given by the canonical refinement of $V_x$.
\end{conjecture}

The case $d=1$ is trivial and
the results above show that this conjectures hold except perhaps for some 
technical restrictions for $d=2$.

The infinitesimal variant of this conjecture (where we replace "there is an affinoid subdomain $\sp S'$ " by "for $S'$ the completion of $S$ at $x$") is a theorem (again under a few technical restrictions): \cite[Theorem 2.5.6]{BC}. The non-criticality assertion is critical in the proof (based on repeated application of Kisin's theorem on $\Lambda^i V$ for $i=1,\dots,d$) and is necessary for the result: see \cite[Remark 2.5.8]{BC}.

Using a similar method (and some argument of Fittings ideals) it is not very hard to prove that after shrinking $X$, we have a filtration $\Fil^i D^\dag(V)$ by sub-families of rank one such that $\Fil^i/\Fil^{i-1}$ are families of rank one. We say that $V$ is {\it weakly trianguline}  If the answer to Question~\ref{Q1} was known to be affirmative, then those families would be {\bf free} of rank one after shrinking $X$,
and after shrinking again (or not, depending on the answer to 
question~\ref{Q2}), they would be 
{\bf of character type} by Lemma~\ref{rankone}, so $D$ would be trianguline.

The proof of the weak triangulinity is omitted here, 
since this result shall be of no use until Question 1 is settled.
Moreover, it seems that Ruochuan Liu has a better proof of that result in preparation, which does not use Kisin's result nor the ugly trick of taking 
exterior powers, but rather works all the way with families of $\fg$-modules and reprove Kisin's theorem as a consequence. It is likely that his proofs shall allow to remove the technical restriction of Theorem~\ref{fam2}. 
It seems also that he could prove Conjecture~\ref{conj} 
up to replacing $X$ by a blow-up. Such a result would be very useful 
for applying Theorem~\ref{main}.

\subsection{Finiteness results for the cohomology of $\fg$-modules}

\label{fin2}

\subsubsection{Cohomology of a general $\fg$-module}
\label{cohfg}
Let $D$ be a general $\fg$-module over $\RR$ (cf.~\ref{different}). 

As in \cite[\S2.1]{Co1}, we define three $\Q_p$-vector spaces $
H^0(D),H^1(D),H^2(D)$ as the cohomology spaces of the complex 
\begin{eqnarray} \label{complex} D \stackrel{d_1}\rightarrow D^2 \stackrel{d2}
\rightarrow D  \end{eqnarray}
with 
\begin{eqnarray}
\label{d1}  d_1(c)&=&((\gamma-1)c,(\phi-1)c), \\
\label{d2} d_2(a,b)&=&(\phi-1)a-(\gamma-1)b.\\
\end{eqnarray}
Of course, the construction of $H^\bullet(D)$ is functorial in $D$, so if $D$ 
has a structure of $S$-module (where $S$ is a ring) compatible with the 
action of $\phi$ and $\Gamma$, then the spaces $H^\bullet(D)$ are naturally $S$-modules as well. Also, even in this generality, it is clear that any short exact sequence of $\fg$-modules $0 \rightarrow D_1 \rightarrow D \rightarrow D_2 \rightarrow 0$ gives rise to a long exact sequence
$$0 \rightarrow H^0(D_1) \rightarrow \dots \rightarrow H^2(D_2) \rightarrow 0.$$
In other words $(H^i)$ is a cohomological functor from the category of 
general $\fg$-modules over $\Q_p$ (resp. over $S$) 
to the category of $\Q_p$-vector spaces (resp. of $S$-modules).
 
\subsubsection{Finiteness results}

\begin{prop}\label{fi}
Let $D$ be a free family of $\fg$-module of rank one of character type over $S$. Then $H^0(D)$ and $H^1(D)$ are finite over $S$.
\end{prop}

\begin{pf} 
We have by hypothesis $D=\RR_S(\delta)$ where $\delta : G_{\Q_p} \rightarrow S^\ast$ is a continuous character. Set $\alpha = \delta(p) \in S^\ast$. The
open sets $\spec S_{a,b}$ defined by the inequalities $a \leq v_p(\alpha(x)) \leq b$ 
for $a,b \in \Z$ form an admissible 
affinoid covering of $S$. Since the functor $\hotimes_S S_{a,b}$ is exact, we get easily, using (\ref{DS'}) that
$$H^i(D) \hotimes_S S_{a,b} = H^i(D_{S_{a,b}})$$ for $i=0,1$.
Therefore, by Kiehl's theorem, it suffices to prove that
$H^i(D_{S_{a,b}})$ are finite over $S_{a,b}$ for $i=0,1$. 
In other words, we  may and shall assume that $a \leq  v_p(\alpha) \leq b$ on $\sp S$ for some integers $a$ and $b$.
 
The proof we shall give is directly inspired by Colmez' computation of the 
$H^0$ and $H^1$ of a $\fg$-module of rank one over a field in \cite[\S2.2]{Co1}.
In his proof, Colmez use several lemmas over the ring $\RR_L$ and $\EE_L^\dag$ 
which, though they are stated only over a field $L$, work as well over a reduced affinoid algebra $S$ since their proofs amount to checking that some explicit series are convergent, which in turn amounts to computing the limits of the norm of their coefficients, a computation that takes exactly the same form in $L$ or in $S$ since in both cases, this
norm is multiplicative, that is, a valuation. This applies to Lemma A.1, A.2, A.4, as well as Corollary A.3 and Lemma 2.3.1 of \cite{Co1}. We shall use those results for $\RR_S$ and $\EE^\dag_S$ without further comments below. 

\par \medskip
{\bf Case of $H^0$\ :}

Since $H^0(\RR_S(\delta)) \subset \RR_S(\delta)^{\phi-1}$, 
and $S$ is Noetherian, the result will follow if we prove that $\RR_S(\delta)^{\phi-1} = \RR_S^{\alpha \phi -1}$ is finite over $S$. 

Actually, we have $\RR_S^{\alpha \phi-1} = (\RR_S^+)^{\alpha \phi-1}$ 
where $\RR_S^{+} = \{\sum_{i \geq 0} a_i T^i \in \RR_S\}$. Indeed, if 
$f = \sum_{i \in \Z} a_i T^i \in \RR_S$ is not in $\RR_S^+$, and $j$ is the largest negative integer $i$ such that $a_i \neq 0$, then one sees, using the fact 
that $\phi(T^{-1}) = T^{-p} + $ terms of lower degrees, 
that if $\phi(f)=\sum b_i T^i$, we have $b_i=0$ for all $i$ such that 
$p j < i < 0$, and therefore one cannot have $\alpha \phi(f)= f$. 

Now choose a non-negative integer $k$ which is strictly greater 
than $-v_p(\alpha)$. The operator $-\sum_{n=0}^{+\infty} (\alpha \phi)^n$ is a continuous inverse of $\alpha \phi-1$ on $T^k \RR_S^+$ by the proof of \cite[Lemma A1]{Co1}.
Therefore the $S$-linear map 
$(\RR_S^{+})^{\alpha \phi=1} \rightarrow \RR_S^+/T^k \RR_S^+$
is injective, and since
the latter has finite rank $k$ over $S$ which is Noetherian, the finiteness of $H^0$ follows.

\par \medskip

{\bf Case of $H^1$, when $v_p(\alpha) < 0$\ :} 

We mimic the beginning of the proof of \cite[Prop. 2.6]{Co1}:

Let $(a,b) \in D^2$, such that $d_2(a,b)=0$ (see~(\ref{d2})). According to \cite[Cor A.3, Lemma A.4 and Lemma 2.3(i)]{Co1}, since $v_p(\alpha)<0$, there exists $c \in D$ such that
if $(a_1,b_1)=(a,b)+d^1(c)$, then $a_1 \in \RR_S^+$, $b_1 \in (\EE_S^\dag)^{\psi=0}$ and
$(\delta (\gamma) \gamma - 1) b_1 \in (\EE_S^+)^{\psi=0}$. In particular, in $D^2$ we have
\begin{eqnarray} \label{D21} \ker d_2 = (\RR_S^+ \times D' ) + \text{Im\,} d^1,\end{eqnarray}
where $D'=\{b \in (\EE_S^\dag)^{\psi=0},
(\delta (\gamma) \gamma - 1) b \in (\EE_S^+)^{\psi=0}\}$. 

Let $k > -v_p(\alpha)$ be an integer.
We claim that in $D^2$ we have the inclusion 
\begin{eqnarray} \label{D22} (T^k \RR_S^+ \times T^k(\EE_S^+)^{\psi=0}) \cap \ker d_2  \subset d_1(\RR_S^+).\end{eqnarray}
Indeed, if $(a,b) \in T^k \RR_S^+ \times T^k (\EE^+)^{\psi=0}$ and $d_2(a,b)=0$,
then by \cite[Lemma A.1]{Co1}, $(\delta(p)\phi-1)$ is an automorphism of $T^k \RR_S^+$, so there exists $c \in \RR_S^+$ such that 
$(\delta(p) \phi -1) c = b$, and since $d_2(a,b)=0$, we have $(\delta(p) \phi-1)((\delta(\gamma)\gamma-1) c -a)=0$, so $(\delta(\gamma)\gamma-1)c=a$ and $d_1(c)=(a,b)$.

Using (\ref{D21}) and (\ref{D22}) we see that there is a surjective $S$-map  
$$(\RR_S^+ \times D')/(T^k \RR_S^+ \times T^k(\EE_S^+)^{\psi=0}) \rightarrow
\ker d_2/\text{Im\, } d_1 = H^1(D).$$ We are thus reduced to prove that $\RR_S^+/T^k \RR_S^+$ is finite over $S$, which is clear, and that $D'/ T^k(\EE_S^+)^{\psi=0}$ is finite. Since $\EE^+_S / T^k \EE^+_S$ is finite, we only need to show that $D'/(\EE^+_S)^{\psi=0}$ is finite, that is, using the definition of $D'$, that the endomorphism $\delta(\gamma)\gamma-1$ of $(\EE^+_S)^{\psi=0}$ has finite cokernel. But this follows from \cite[Lemma 2.3(ii)]{Co1}.

\medskip

{\bf Case of $H^1$, general $v_p(\alpha)$ :} 

As in \cite[\S 2.4]{Co1}, there is an $S$-linear operator $\partial : H^1 (\RR_S(t^{-1}\delta)) \rightarrow H^1(\RR_S(\delta))$ where $t$ is the character $\Q_p^\ast \rightarrow S^\ast, \ t \mapsto t$. Mimicking \cite[S 2.4]{Co1}, it is easy to see that $\coker\,
\partial$ is finite over $S$. Therefore, $H^1(\RR_S(\delta))$ is finite if $H^1(\RR_S(t^{-1} \delta))$ is. By induction, $H^1(\RR_S(\delta))$ is finite if $H^1(\RR_S(t^{-n} \delta))$ is
for some $n \in \N$. But since $v_p( \delta(p) ) \geq b$ on $S$, taking $n \geq b+1$, we get that $H^1(\RR_S(t^{-n} \delta))$ by the preceding case, so $H^1(\RR_S(\delta))$ is finite.
\end{pf}

\begin{remark} It is not hard, actually even simpler than for the $H^1$,
to prove that $H^2(D)$ is finite over $S$ as well. We leave this as an exercise for the reader.
\end{remark}

\begin{cor}
\label{fincor2}Let $D$ be a trianguline family of $\fg$-module over $S$.
Then $H^0(D)$ and $H^1(D)$ are finite over $S$.
\end{cor}
\begin{pf} This is obvious by the long exact sequence of cohomology and the
preceding proposition.
\end{pf}

\subsection{The local Pottharst Selmer group}

\label{localpot}

Let $S$ be a reduced affinoid algebra over $\Q_p$, $V$ a trianguline family of representation over $S$ of rank $d$, and $a$ an integer in $\{0,\dots,d\}$.
Following Pottharst (\cite{pot}) we define an $S$-linear application 
\begin{eqnarray*} u : H^1(G_{\Q_p},V) = \Ext^1_{G_{\Q_p}}(1,V) \stackrel{D^\dag}{\rightarrow} 
\Ext^1_{\fg-\text{modules}}(1,D^\dag(V))\\ = H^1(D^\dag(V)) \rightarrow H^1(D^\dag(V)/
\Fil^a D^\dag(V).\end{eqnarray*}

 Above, $1$ denotes both the trivial representation of dimension $1$ over $S$
and the trivial $(\phi,\Gamma)$-module over $\RR_S$;  the morphism denoted 
$D^\dag$ is the morphism induced by the functor $D^\dag$ using the fact that $D^\dag(1)=1$ and that $D^\dag$ transforms exact sequences into exacts sequences, 
therefore extension of $1$ by $V$ into extensions of $1$ by $D^\dag(V)$; 
the last map is the $H^1$ of the canonical surjection.
We define $H^1_{\pot,a}(G_{\Q_p},V) = \ker (H^1(G_{\Q_p},V) \stackrel{u} 
\rightarrow H^1(D^\dag(V)/\Fil^a D^\dag(V)))$ 

The following result is a close variant of the main theorem of Pottharst' paper 
\cite[Theorem 3.1]{pot}. The hypotheses and the proof we give is different, relying
only on results of \cite[Chapter II]{BC}, and close in spirit to the  proof of Flach's theorem relating Greenberg and Bloch-Kato Selmer groups. 

\begin{prop} \label{pot}
Assume that $S$ is a finite extension $L$ of $\Q_p$,  
that $V$ is a non-critical and non-exceptional trianguline representation over 
$L$, and that $a$ is the number of non-positive 
Hodge-Tate weights of $V$. Then $H^1_{\pot,a}(G_{\Q_p},V)=H^1_{g}(G_{\Q_p,V})=H^1_f(G_{\Q_p},V)=H^1_e(G_{\Q_p},V)$.
\end{prop}
\begin{pf}
We first show that $H^1_g=H^1_f=H^1_e$ for a trianguline non-exceptional representation. By \cite[Corollary 3.8.4]{BK} (for $H^1_e = H^1_f$) and 
\cite[Proposition 3.8]{BK} (for $H^1_g=H^1_e$), 
it is enough to show that $1$ and $p$ are  not eigenvalues of $\phi$ on $D_\crys(V)$. Let $(\delta_i)_{i=1,\dots,d}$ be the parameter of $V$ and let us note 
$\D_\crys$ the functor from the category of $\fg$-modules over $L$ to the category of filtered $\phi$-modules (such that $D_\crys = \D_\crys \circ D^\dag$)
constructed by Berger, namely $\D_\crys(D)=D[1/t]^\Gamma$. By an immediate d\'evissage, it is enough to show that $1$ and $p$ are not eigenvalue
of $\phi$ 
on $\D^\crys(\delta_i)$ (a space of dimension $1$ or $0$, in which case it has no $\phi$-eigenvalue at all).
But we see at once that this space is not $0$ if and only if $\delta_i(\gamma)=\gamma^{-n}$ for $\gamma \in \Gamma$, where $n$ is a rational integer, 
necessarily equal to the weight $s_i$ of the character $\delta_i$. In this case, the action of $\phi$ is $p^{n} \delta(p)$, so we see that the only case where this can be one or $p$ is when $\delta_i(t) = t^n$ or $\delta_i(t)=t^{n}|t|$ for all $t \in \Q_p^\ast$, that is when $\delta_i$ is exceptional.

We now prove that $H^1_{\pot,a}(V) \subset H^1_{g}(G_{\Q_p},V)$ for a trianguline non critical $V$. We have to show that an extension $U$ (in the category of $\fg$-module) of $1$ 
by $D^\dag(V)$ whose push-forward as an extension of $1$ by 
$D^\dag(V)/\Fil^a D^\dag(V)$ is trivial, is de Rham. But this amounts to 
showing that an extension of
$1$ by $\Fil^a D^\dag(V)$ is de Rham (as a $\fg$-module). 
Such an extension is a trianguline $\fg$-module of rank $a+1$ which has for 
parameter 
the sequence of characters $\delta_1,\dots,\delta_a,1$, 
whose weights are $s_1,\dots,s_a,0$. 
Since $V$ is non-critical, and since $a$ is its number of 
non-positive Hodge-Tate weights, we have $s_1<s_2<\dots<s_a<0$, 
and this $\fg$-module is de Rham by \cite[Proposition 2.3.4]{BC}.

Finally, we conclude the proof by showing  that 
$\dim H^1_{\pot,a}(G_{\Q_p},V)=\dim H^1_f(G_{\Q_p},V)$ for a trianguline, non-exceptional $V$. 
Indeed, since $D^\dag$ is an equivalence of category (we are over a field $L$),
it realizes an isomorphism $H^1(G_{\Q_p},V) \simeq H^1(D^\dag(V))$ which identifies $ H^1_{\pot,a}(G_{\Q_p},V)$ with $\ker (H^1(D^\dag(V)) \rightarrow H^1(D^\dag(V)/Fil^a D^\dag(V)))$.
Since the characters $\delta_i$ in the parameter of $V$ are non-exceptional, we have $\dim H^i(\RR_L(\delta_i)) =0$ if $i=0$ or $2$
and $1$ if $i=1$. In particular, we get easily by d\'evissage (and descending induction on $i$, from $i=d$ to $i=0$) that
 $H^0( D^\dag(V)/\Fil^i D^\dag(V))=0$. Applying this for $i=a$ 
we get using the long exact sequence attached to $0 \rightarrow \Fil^a D^\dag(V) \rightarrow D^\dag(V) \rightarrow D^\dag(V)/\Fil^a D^\dag(V) \rightarrow 0$ 
that $$H^1_{\pot,a}(G_{\Q_p}, V) \simeq H^1(\Fil^a D^{\dag}(V))$$ 
and a new d\'evissage shows that this space has dimension $a$. On the other hand, since as we have seen $H^0(D^\dag(V))=0$, so $H^0(G_{\Q_p},V)=0$, the dimension of $H^1_f(G_{\Q_p},V)$ is by \cite[Corollary 3.8.4]{BK} the number of non-positive Hodge-tate weights of $V$, that is, by hypothesis, $a$. 

\end{pf}

\section{Proof of Theorem~\ref{main}}

Let $K$ be a number field, $p$ a prime number that is totally split in $K$, 
$\Sigma_p$ the set of places of $K$ 
above $p$, $\Sigma$ a finite set of places containing $\Sigma_p$.

For any affinoid algebra $S$ over $\Q_p$, and any locally free module $V$ of rank $d$ over $S$ with a continuous action of $G_K$, such that 
$D_v := D^\dag(V_{|G_v})$ is trianguline (that is, 
provided with a triangulation $\Fil^\cdot D_v$) and any
fixed sequences of integers $\aa = (a_v)_{v \in \Sigma_p}$, $0 \leq a_v \leq d$ we define the global {\it Pottharst Selmer group}
as \begin{eqnarray*} 
H^1_{\pot,\aa} (G_{K,\Sigma},V) &=& \ker (H^1(G_{K,\Sigma}, V) \stackrel{u}\rightarrow \prod_{v \in \Sigma_v} H^1(G_v, D_v/\Fil^{a_v}(D_v))\\
&=& \{x \in H^1(G_{K,\Sigma},V),\ \forall v \in \Sigma_p,\  x_v \in H^1_{\pot,a_v}(G_v,V)\}\end{eqnarray*}
In the first line, $u$ is the composition of restriction map 
$$H^1(G_{K,\Sigma},V) \rightarrow \prod_{v \in \Sigma_p} H^1(G_v,V_{|G_v})$$ with the
 products of the maps $u_v : H^1(G_v,V_{|G_v}) \rightarrow H^1(D_v/Fil^{a_v}D_v)$ constructed (and denoted $u$) in \S\ref{localpot}.
This notion is taken from \cite{pot}, where it is called 
{\it triangulordinary Selmer group}

Over a field, we have the following interpretation:
\begin{prop} \label{globalpot}
Assume that $S=L$ is a finite extension of $\Q_p$, and assume that for each $v \in \Sigma_p$,
$V_{|G_v}$ is trianguline, non critical, and non exceptional, and that 
$a_v$ is the number of non-positive Hodge-Tate weights of $V_{|G_v}$ Then 
$H^1_{\pot,\aa}(G_{K,\Sigma},V)=H^1_g(G_{K,\Sigma},V)=H^1_f(G_{K,\Sigma},V)=H^1_e(G_{K,\Sigma},V)$
\end{prop}
\begin{pf} 
This follows immediately from the definitions and Prop.~\ref{pot}.
\end{pf}
 
We want to prove the main theorem. So assume that $S$ is a domain of 
dimension $1$. It is clearly allowed to replace $S$ by its normalization, 
so we assume that $S$ is normal.
Let $x$ be a point of $\sp S$. It corresponds to a prime principal ideal $(f)$.
Then $S/f=L(x)$, the field of definition of $x$, and applying the 
construction above to $S/f$ and $V/fV$, we get a map 
$u_f : H^1(G_{K,\Sigma},V/fV) \rightarrow \prod H^1(D^v / Fil^a D_v)$. We are now in position to apply the formalism developed 
in~\ref{formalism}. Indeed we can take for $\Ac$ 
the category of finite $S$-modules endowed with a continuous representation
 of $G_{K_\Sigma}$, and for $(H^i_\Ac)_{\i \in \N}$, we take the usual continuous cohomology groups, exactly as in \S\ref{catAc}. As for $\Ac'$, we take
 the product, indexed by $\Sigma_p$,  of copies
the category of general $(\phi,\Gamma)$-modules (cf. \S\ref{different}), and for $H^i_{\Ac'}$ we take the product over $\Sigma_p$ of the cohomology of 
$\fg$-modules (cf.~\S\ref{cohfg})  
 
For $V$ we take our $S$-representation $V$, and for $W'$ we take 
$(D_v/\Fil^{a_v} D_v)_{v \in \Sigma_p} \in \Ac'$.
The maps $u : H^1(V) \rightarrow H^1(V')$ and 
$u_f : H^1_f(V/fV) \rightarrow H^1(V'/fV')$ have been defined above, and the commutativity of the 
diagram~(\ref{commute}) is trivial if we take into account the fact that $D^\dag(V/fV)$ is naturally isomorphic to $D^\dag(V)/fD^\dag(V)$ by Berger-Colmez theorem.
By definition, $S(V)=\ker u = H^1_{\pot,\aa}(G_{K,\Sigma},V)$ and $S(V/fV)=\ker u_f = H^1_{\pot,\aa}(G_{K,\Sigma},V/fV)$.

We want to apply Corollary~\ref{form3} to this situation. We only have to check the hypotheses that $S^1(G_{K,S},V)$, $H^2(G_{K,S},V)$ and $\prod_{v \in \Sigma_v} 
H^1(D_v/\Fil^{a_v} D_v)$ are finite over $S$. But the two first follow from Prop~\ref{finiteness} while the second is Corollary~\ref{fincor2}.
Corollary~\ref{form3} thus tells us that setting $r = r_{(a_v)} = \rk_S S(V)$ we have for all prime principal ideal $(f)$ of $S$:  
$$\dim S(V/fV) \geq r,$$ withe equality for all $(f)$ except a finite number. The same of course will holds if we restrict our attention
to the subset of  $(f)$'s such that $V/fV$ is, at each $v \in \Sigma_K$, 
non-critical, non-exceptional, and with exactly $a_v$ non-positive Hodge-Tate weights. But for those $(f)$, we know by the Proposition~\ref{globalpot} 
that $S(V/fV)=H^1_{\pot,\aa}(V/fV)=H^1_g(G_{K,\Sigma},V/fV)= H^1_f(G_{K,\Sigma},V/fV)= 
H^1_e(G_{K,\Sigma},V/fV)$ and this complete the proof of the theorem.

\begin{remark} If we had a positive answer to Question~\ref{QDdag} (resp. to Question 
~\ref{QDdag2}), we could apply (resp. adapt) the formalism of \S\ref{formalism2} with $\Ac$
the category of finite $S$-module with $G_{K,\Sigma}$-action, and $\Ac'$ the category of general $\fg$-modules over $S$. This would prove that the statement of 
Theorem~\ref{main}
would hold for the regular points of any reduced affinoid $S$, not necessarily one of dimension $1$.
\end{remark}
\section{Application to modular forms}

In \cite{BC} the following result, a special case of the Bloch-Kato conjecture, is proved.\footnote{Actually, it is proved under two assumptions, called there 
Rep(4) and AC$(\pi_{f,E})$ which at that time were 
announced in the literature, but a complete proof of which had not appeared yet. The situation is a little bit better now, as preprints of  Morel, Shin, and a team led by Harris cover what it needs to give a proof of Rep$(4)$ (and Rep$(n)$ for that matter). cf. \cite{morel}, \cite{shin} and \cite{book}} 

\begin{theorem}\cite[Cor. 8.1.4]{BC} \label{signe}
Let $f = \sum_{n \geq 1} a_n q^n$ be a newform of weight $k=2k' \geq 4$ and level $\Gamma_0(N)$, with $N$ prime to $p$.
Let $\rho=\rho_f(k')$ be the Galois representation of $G_{\Q,\Sigma}$ attached to $f$
 (where $\Sigma$ is the set of primes dividng $Np$) normalized so to have motivic weight $-1$.
Then if $\varepsilon(\rho,0)=-1$ we have  $\dim H^1_f(G_{\Q,\Sigma},\rho) \geq 1$.
\end{theorem}
We want to extend this result:
\begin{theorem} \label{signe2} The theorem above also holds if $k=2$ except maybe if $a_p=2\sqrt{p}$.

 In particular, if $E$ is an elliptic curve over $\Q$ with good reduction at $p$,
and such that $\epsilon(E,0)=-1$, then the $p$-adic Selmer group of $E$ has 
rank at least one.
\end{theorem}
Of course, the ``in particular'' follows from the first sentence using the fact that for our elliptic curve $E$, there is a modular newform $f$
of weight $2$ with integral coefficients (in particular, with 
$a_p \neq 2\sqrt{p}$) and level $\Gamma_0(N)$ such that $\rho_{f}(1) = V_p(E)$ and the fact that the Bloch-Kato Selmer group
$H^1_f(\Q,V_p(E))$ is the same as $\Sel_p(E) \otimes_{\Z_p} \Q_p$ where 
$\Sel_p(E)$ is the $p$-adic Selmer group of $E$.

\begin{remark} 
\begin{itemize} 
\item[(i)]
In the case where $p$ is an ordinary prime for $f$, there already exists two other proof of this 
theorem (and of Theorem~\ref{signe2}): one due to Nekovar and one due to Skinner and Urban. 

\item[(ii)] In the case of an elliptic curve $E$, there is another proof of that result, by completely different method, due to Kim. See \cite{kim}.

\item[(iii)] The technical hypothesis $a_p \neq 2\sqrt{p}$ arises from
the technical hypothesis in our triangulinity result Theorem~\ref{fam2}. 
A better version of that theorem should allow to remove that hypothesis.
\end{itemize}
\end{remark}

Let us prove Theorem~\ref{signe2}. In view of Nekovar and Skinner-Urban's result, we may assume that $f$ is supersingular at $p$. (Actually, we could treat the ordinary case with the same method using Remark~\ref{remfam2} instead of theorem ~\ref{fam2}), but that would require a boring case-by-case analysis).

The representation $(\rho_f(1))_{|G_p}$, which is crystalline and irreducible, 
always admits two non-critical refinements $\RR$ (\cite[Remark 2.4.6(iii)]{BC}). 
Let us choose one of them, so we fix an ordering of the crystalline 
eigenvalues $(\phi_1,\phi_2)$. The $\fg$-module $D=D^\rig((\rho_f(1))_{|G_p})$ is thus trianguline, of parameter $(\delta_1,\delta_2)$ with $\delta_1(p)=p \phi_i$, $\delta_2(p)=\phi_2$ and 
$\delta_1(t)=t$, $\delta_2(t)=1$ for $t \in \Gamma=\Z_p^\ast$ (using the fact that the Hodge-Tate weights of $\rho_f(1)_{|G_{\Q_p}}$ are $-1$ and $0$ -- see 
\cite[Prop. 2.4.1]{BC}).
By weak admissibility, and irreducibility, we have $-1 < v_p(\phi_i) < 0$ for 
$i=1,2$. Therefore, $\delta_1$ and $\delta_2$, hence $\rho_f(1)_{|\G_{\Q_p}}$, 
are non-exceptional. Moreover, we see that the only case 
where $\delta_2 \delta_1^{-1} (p) = \phi_2/(p \phi_1)$ may be in $p^\Z$ are when $\phi_2=\phi_1$. Since the product $\phi_1 \phi_2$ is $p^{-1}$ and the sum $\phi_1+\phi_2$ is $a_p p^{-1}$, this is equivalent to $a_p = 2 \sqrt{p}$, so does not happen under our hypotheses.

Let $\EE$ denotes the eigencurve of tame level $N$. Then $(f,\RR)$
is a point $x$ of the $\EE$. Let $X$ be an irreducible affinoid neighborhood of $x$ in an irreducible component of $\EE$ through $x$. Then $X=\sp S$ where $S$ is an affinoid domain of dimension $1$. 

Since $\rho_f$ is irreducible, by Rouquier's theorem, there exists (after replacing the base field $\Q_p$ by a finite extension $L$ if necessary, and shrinking $X$) 
 a Galois representation $V$ of  $G_{\Q,\Sigma}$ ($\Sigma$ being the set of rational primes dividing $Np$) of dimension $2$ on $X$
whose trace is the the restriction to $X$ of the canonical pseudocharacter carried by $\EE$,
normalized so that at a point $z \in X$ corresponding to a {\bf classical} 
modular eigenform $g_z$ of weight $k_z=2k'_z$ and of level $\Gamma_0(N)$, we have $V_z \simeq\rho_{g_z}(k'_z)$. We call $Z$ the set of such classical points. The family 
$V_{|G_{\Q_p}}$  over $X$ together with the set of classical points $Z$ is
a refined family in the sense of~\ref{refined}. This follows from the construction of the eigencurve, see~\cite[Prop. 7.5.13]{BC} for details. Note that
for $z \in Z$ attached to a form $g_z$ of weight $2k'_z$, the Hodge-Tate weights of 
$V_z$ are $-k'_z$ and $k'_z-1$, so $V_z$ has exactly $1$ non-positive weight.

By Theorem~\ref{fam2} (applicable since $(V_x)_{|G_{\Q_p}}$ is irreducible, trianguline and non-critical, and
its parameter satisfies $(\delta_2 \delta_1^{-1})(p) \not \in p^\Z$), up to shrinking $\sp S$ if necessary, the family $V_{|G_p}$ is trianguline 
(the triangulation of $(V_x)_{|G_p} = (\rho_f)_{|G_p}$ being associated with the chosen refinement). We now apply 
Theorem~\ref{main}, for $\aa=(a_v)_{v \in \Sigma_p=\{p\}}=(1)$: we get that there 
exists an integer $r$ such that for all $x \in \sp S$ such that $(V_x)_{|G_p}$
is non-critical, non-exceptional, and with exactly one non-negative Hodge-Tate weight, $\dim H^1(G,V_x) \geq r$ with equality 
 up to a finite number of exception. We apply this to the set of points $Z' \subset Z$ of $z$ that correspond to classical forms $g$ of weight $2k' \geq 4$ 
with a non-critical, non-exceptional, refinement. 
It is easily seen using (v) of the definition of a refined family (\S\ref{refined}) that $Z'$ is infinite. Therefore for an infinity of $z \in Z'$, 
we have $\dim H^1_f(G_{\Q_,\Sigma},V_z)=r$.
But as is well known, the sign $\varepsilon(0)$ is constant in a connected
family, so by hypothesis, if a form $g$ of weight $2k'$ 
corresponds to a $z \in Z'$, with $\rho_g(k') \simeq V_z$,
 we have $\varepsilon(\rho_g(k'),0)=-1$ and
by Theorem~\ref{signe2}, $\dim H^1_f(G_{\Q,\Sigma},\rho_g(k')) \geq 1$ 
for those $g$, so $r \geq 1$. It follows that $\dim H^1_f(G_{\Q,\Sigma},\rho_f(1)) \geq 1$.

\par \bigskip
 
 { }

\end{document}